\documentclass[a4paper,12pt]{article}
\usepackage{amsbsy,amsfonts,amsmath,amssymb,amsthm}
\usepackage{array}
\usepackage{bm}
\usepackage{color}
\usepackage{enumerate}
\usepackage{mathrsfs}
\usepackage{latexsym}
\usepackage[hidelinks]{hyperref}
\usepackage{mathtools}

\mathtoolsset{showonlyrefs,showmanualtags}

\definecolor{wineRed}{rgb}{0.7,0,0.3}
\definecolor{granduncle}{rgb}{0,0,0.8}
\definecolor{Green}{rgb}{0,0.4,0}
\definecolor{blueViolet}{rgb}{0.4,0,1.0}
\definecolor{bloodOrange}{rgb}{0.85,0.05,0}
\definecolor{mycolor}{rgb}{0.8,0,0.2}
\definecolor{}{rgb}{0.8,0,0.2}

\usepackage[textsize=small]{todonotes}
\usepackage{cite}

\usepackage{comment}
 
\DeclareMathAlphabet{\mathpzc}{OT1}{pzc}{m}{it}

\numberwithin{equation}{section}

\theoremstyle{plain}
 
\newtheorem{lem}{Lemma}

\newtheorem{prop}{Proposition}

\theoremstyle{definition}
\newtheorem{defn}{Definition}
\newtheorem{thm}{Theorem}
\newtheorem{rem}{Remark}

\newtheorem{ex}{Example}

\def\N{\mathbb{N}}
\def\R{\mathbb{R}}

\def\ds{\displaystyle}
\def\ts{\textstyle}

\begin{document}
                   \begin{center}
{\Large \bf Periodic Solutions to \\[3mm]
 Kobayashi--Warren--Carter Systems}\\[1cm]
                   \end{center}
               \vspace{0.7cm}

                  \begin{center}
   \textsc{Shodai Kubota}\\
Department of Mathematics, Faculty of Engineering, Kanagawa University\\
3-27-1, Rokkakubashi, Kanagawa-ku, Yokohama, 221-8686, Japan\\
E-mail : skubota@kanagawa-u.ac.jp
                   \vspace{0.7cm}

                  {\sc Ken Shirakawa}\\
Department of Mathematics, Faculty of Education, Chiba University\\
1-33 Yayoi-cho, Inage-ku, Chiba, 263-8522, Japan\\
E-mail: sirakawa@faculty.chiba-u.jp
                   \vspace{0.7cm}
\end{center}

\hspace*{-0.6cm}{\bf Abstract.} 
In this paper, a system of parabolic PDEs, called the Kobayashi--Warren--Carter system, is considered as a possible phase-field model of planar grain boundary motion. The Main Theorem is concerned with the existence of a time-periodic solution to the Kobayashi--Warren--Carter system, and the principal objective is to provide a proof without the use of a compromised assumption, which researchers have been forced to adopt in recent studies.  
\\[0.2051cm]

\hspace*{-0.6cm}{\bf AMS Subject Classification:} 
35B10,  
35K59,  
35K67,  
35K87.  
\\[0.2051cm]

\hspace*{-0.6cm}{\bf Keywords:} 
Kobayashi--Warren--Carter system, grain boundary motion, time-periodic solution.\\

\newpage

\section{Introduction}
Let $ (0, T) $ be a time interval with a constant $ 0 < T < \infty $ and let $N \in \mathbb{N}$ denote the spatial dimension. 
Let $\Omega \subset \R^N $ be a bounded domain with a Lipschitz boundary $ \Gamma := \partial \Omega $ and let $ n_\Gamma $ be the unit outer normal on $ \Gamma $. Additionally, we set $ Q := (0, T) \times \Omega $ and $ \Sigma := (0, T) \times \Gamma $. 

In this paper, we consider the following system of parabolic PDEs, denoted by (S):

~~~~\hypertarget{(S)}{(S)}
\begin{equation}\label{1}
\left\{ \parbox{11cm}{
    $ \partial_{t} \eta - \kappa^2 \mathit{\Delta} \eta +g(\eta) +\alpha'(\eta) |D \theta| =  u  $ \quad in $ Q $,
\\[1ex]
$ \nabla \eta \cdot n_\Gamma = 0 $ \quad on $ \Sigma $, 
}\right. 
\end{equation}
\begin{equation}\label{2}
\left\{\parbox{11cm}{
    $ \displaystyle \alpha_{0}(\eta) \partial_{t} \theta -\mbox{div} \left( \alpha(\eta) \frac{D \theta}{|D \theta|} +\nu_0 \nabla (\nu_0 \theta) \right) + M_0\theta =  v  $ \quad in $Q$,
\\[1ex]
    $\displaystyle  \bigl( {\ts \alpha(\eta) \frac{D \theta}{|D \theta|} +\nu_0 \nabla (\nu_0 \theta)} \bigr) \cdot n_\Gamma = 0 $ \quad on $ \Sigma $. 
}\right. 
\end{equation}
This system is based on the Kobayashi--Warren--Carter system, which was proposed by Kobayashi--Warren--Carter \cite{MR1752970,MR1794359}, as a possible phase-field model of planar grain boundary motion. According to the original studies \cite{MR1752970,MR1794359}, it is supposed that the system (S) is governed by the following \emph{free energy,} denoted by $ \mathcal{F}$:
\begin{align}\label{freeEnergyOrg}
    \mathcal{F} : [{\eta}, \,& {\theta}] \in H^1(\Omega) \times BV(\Omega) \mapsto  \mathcal{F}({\eta}, {\theta}) := \frac{\kappa^2}{2} \int_\Omega |\nabla {\eta}|^2 \, dx +\int_\Omega G({\eta}) \, dx
    \nonumber
    \\
    & +\int_\Omega \alpha({\eta}) |D {\theta}| +\frac{1}{2} \int_\Omega |\nabla (\nu_0 {\theta})|^2 \, dx +\frac{M_0}{2} \int_\Omega |{\theta}|^2 \, dx \in [0, \infty],
\end{align}
and the polycrystalline microstructure, including the grain boundary, is reproduced by a vector field $ \eta \bigl[ \cos \theta, \sin \theta \bigr] $ in $ Q $, which consists of two unknown variables $ \eta $ and $ \theta $. In this context, $ \eta = \eta(t, x) $ and $ \theta = \theta(t, x) $ are the order parameters of the orientation order and orientation angle of crystallization, respectively. $ \kappa > 0 $ is a fixed constant of the spatial diffusion of $ \eta $. $ g = g(\eta) $ is a Lipschitz perturbation used to control the range of $ \eta $ and $ G = G(\eta) $ is a nonnegative potential of $ g $. $ \alpha = \alpha(\eta) $ and $ \alpha_0 = \alpha_0(\eta) $ are given positive-valued functions that represent the mobilities of the grain boundary. $ M_0 > 0 $ and $ \nu_0 \geq 0 $ are fixed constants. Finally, $ u = u(t, x) $ and $ v = v(t, x) $ are given forcing terms for $ \eta $ and $ \theta $, respectively, and in particular, $ u $ is supposed to correspond to the relative temperature in $ Q $. 

In the system (S), the principal target is the piecewise constant case of the orientation angle $ \theta $, which is achieved by the weighted total variation $ \int_\Omega \alpha(\eta)|D \theta| $ in \eqref{freeEnergyOrg}, and it is desirable to capture the \emph{facet structure} in a polycrystal. In the original studies by Kobayashi--Warren--Carter \cite{MR1752970,MR1794359}, a number of targeted solutions were exemplified via structural observations in the 1D setting of $ \Omega $, and numerical experiments in the 1D and 2D settings of spatial domains. Additionally, Giga--Giga \cite{MR2746654} and Giga--Kobayashi \cite{MR1865089} addressed the 1D-singular diffusion equation as a toy problem of \eqref{2} and successfully analyzed the fundamental rule of the behavior of piecewise constant solutions.     

Simultaneously , the mathematical analysis of general solutions to (S) (including the targeted solutions) was studied in previous works \cite{MR2836555,MR2469586,MR2548486,MR2668289,MR3268865,MR3038131,MR3203495,MR3238848,MR3670006,MR3082861,MR4218112,KNS2020,MR4228007} using the abstract theory of nonlinear evolution equations. Most of the previous studies were devoted to the existence and large-time behavior of solutions, and are broadly divided into two categories: 
\begin{itemize}
    \item[$\sharp 1)$] results in the standard case when $ \nu_0 > 0 $, which is based on the standard variational method of a nonlinear parabolic system (cf. studies by Ito, Kenmochi, and Yamazaki \cite{MR2836555,MR2469586,MR2548486,MR2668289}); 
    \item[$\sharp 2)$] results in the non-standard case when $ \nu_0 = 0 $, which developed the mathematical treatment of the weighted total variation measure $ \alpha'(\eta)|D \theta| $ in \eqref{1} and the singular diffusion $ -\mathrm{div} \bigl( \alpha(\eta) \frac{D \theta}{|D \theta|} \bigr) $ in \eqref{2} (cf. studies by Moll,  Nakayashiki, Watanabe, Yamazaki, and the authors of the present paper \cite{MR3268865,MR3038131,MR3203495,MR3238848,MR3670006,MR3082861,MR4218112,KNS2020,MR4228007}). 
\end{itemize}
We identified a few results concerned with the uniqueness of (S), and currently, any uniqueness result is obtained only in the case when $ \nu_0 > 0 $. Furthermore, except for the 1D case of $ \Omega $ (cf. \cite[Theorem 2.2]{MR2469586}), we require the following compromised assumption for the uniqueness of the Kobayashi--Warren--Carter system (S):
\begin{itemize}
    \item[$\sharp3)$]the mobility $ \alpha_0 $ is constant with respect to $ \eta $. 
\end{itemize}

In recent studies on the system (S), the researchers have been advanced to the issues of control/stabilization problems of grain boundary motion, and we identified several results for optimization problems (cf. \cite{MR2668289,MR4218112,KNS2020,MR4228007}). However, because of the difficulty of uniqueness, the researchers that obtained these results were forced to rely on the compromised assumption $\sharp3)$; hence, current results are insufficient as a complete mathematical theory in the field of nonlinear science. 

In view of this, we focus on the $T$-periodic solution to (S), that is, the time-periodic solution with period $ T $, as fundamental research on the \emph{stabilization} of grain boundary motion. The Main Theorem is concerned with the existence of a $T$-periodic solution and our principal objective is to prove the Main Theorem in the original $ \eta $-dependent setting of $ \alpha_0 $. 

The structure of this paper is as follows: We state the Main Theorem in Section 3 based on the preliminaries in Section 2. In Section 4, we consider the time-discretization scheme for (S), as the approximating problem, and obtain the time-periodic solution to the approximating problem using a priori estimates. Finally, in Section 5, we prove the Main Theorem as a consequence of the limiting observation of approximating the time-periodic solution. Furthermore, we present some elementary tools in the appendix, such as the Gronwall-type inequality in time discretization.

\section{Preliminaries}

We begin by prescribing the notations used throughout this paper. 
\medskip

\noindent
\underline{\textbf{\textit{Basic notations.}}} 
For arbitrary $ r_0 $, $ s_0 \in [-\infty, \infty]$, we define:
\begin{equation*}
r_0 \vee s_0 := \max\{r_0, s_0 \}\ \mbox{and}\ r_0 \wedge s_0 := \min\{r_0, s_0 \},
\end{equation*}
and in particular, we set:
\begin{equation*}
    [r]^+ := r \vee 0 \ \mbox{and}\ [r]^- :=  -(r \wedge 0), \mbox{ for any $ r \in \R $.}
\end{equation*}    
For any dimension $ k \in \N $ and any set $ A $, we write:
\begin{equation*}
    [A]^k := \overbrace{A \times \cdots \times A}^{\mbox{$k$ times}}.
\end{equation*}

\noindent
\underline{\textbf{\textit{Abstract notations.}}}
For an abstract Banach space $ E $, we denote by $ |\cdot|_{E} $ the norm of $ E $, and denote by $ \langle \cdot, \cdot \rangle_E $ the duality pairing between $ E $ and its dual $ E^* $. 
In particular, when $ H $ is a Hilbert space, we denote by $ (\cdot,\cdot)_{H} $ the inner product of $ H $. 

For a proper functional $ \Psi : E \longrightarrow (-\infty, \infty] $ on a Banach space $ E $, we denote by $ D(\Psi) $ the domain of $ \Psi $, i.e. 
    $ D(\Psi) := \left\{ \begin{array}{l|l}
        z \in E & \Psi(z) < \infty
    \end{array} \right\} $.
\medskip

\noindent
\underline{\textbf{\textit{Notations in variational analysis. (cf. \cite{MR0348562,MR0298508,MR1201152})}}} 
    Let $ H $ be an abstract Hilbert space. For a proper, lower semi-continuous (l.s.c.), and convex function $ \Psi : H \longrightarrow (-\infty, \infty] $, 
    we denote by $\partial \Psi$ the subdifferential of $\Psi$. The subdifferential $ \partial \Psi $ corresponds to a weak differential of convex function $ \Psi $, and it is known as a maximal monotone graph in the product space $ H \times H $. The set $ D(\partial \Psi) := \bigl\{ z \in H \ |\ \partial \Psi(z) \neq \emptyset \bigr\} $ is called the domain of $ \partial \Psi $. We often use the notation ``$ [z_{0}, z_{0}^{*}] \in \partial \Psi $ in $ H \times H $\,'', to mean that ``$ z_{0}^{*} \in \partial \Psi(z_{0})$ in $ H $ for $ z_{0} \in D(\partial\Psi) $'', by identifying the operator $ \partial \Psi $ with its graph in $ H \times H $.
\medskip

\begin{ex}[Examples of the subdifferential]\label{exConvex}
    As one of the representatives of the subdifferentials, we exemplify the following set-valued function $ \bigl[ \frac{(\cdot)}{|\,\cdot\,|} \bigr]: \R^d \longrightarrow 2^{\mathbb{R}^d} $, with $ d \in \mathbb{N} $, which is defined as:
\begin{align}\label{Sgn^d}
    \xi = [\xi_1, & \dots, \xi_d] \in \mathbb{R}^d \mapsto \bigl[ {\ts\frac{\xi}{|\xi|}} \bigr] 
    := \left\{ \begin{array}{ll}
            \multicolumn{2}{l}{
                    \ds \frac{\xi}{|\xi|} = \frac{[\xi_1, \dots, \xi_d]}{\sqrt{\xi_1^2 +\cdots +\xi_d^2}}, ~ } \mbox{if $ \xi \ne 0 $,}
                    \\[3ex]
            \mathbb{D}^d, & \mbox{otherwise,}
        \end{array} \right.
    \end{align}
    where $ \mathbb{D}^d $ denotes the closed unit ball in $ \mathbb{R}^d $ centered at the origin. Indeed, the set-valued function $ \bigl[ {\ts\frac{(\cdot)}{|\,\cdot\,|}} \bigr] $ coincides with the subdifferential of the Euclidean norm $ |{}\cdot{}| : \xi \in \mathbb{R}^d \mapsto |\xi| = \sqrt{\xi_1^2 + \cdots +\xi_d^2} \in [0, \infty) $, i.e.:
\begin{equation*}
    \partial |{}\cdot{}|(\xi) = \bigl[ {\ts\frac{\xi}{|\xi|}} \bigr], \mbox{ for any $ \xi \in \mathbb{R}^d ~\bigl( = D(\partial |{}\cdot{}|) \bigr) $.}
\end{equation*}
\end{ex}
\begin{ex}\label{Rem.gamma_eps}
    Let $ d \in \mathbb{N} $ be the constant of dimension. For any $ \varepsilon \geq 0 $, let $ \gamma_\varepsilon : \mathbb{R}^d \longrightarrow [0, \infty) $ be a continuous and convex function, defined as:
    \begin{equation}\label{gamma_eps}
        \gamma_\varepsilon : y \in \mathbb{R}^d \mapsto \gamma_\varepsilon(y) := \sqrt{\varepsilon^2 +|y|^2} \in [0, \infty).
    \end{equation}

    When $ \varepsilon = 0 $, the convex function $ \gamma_0 $ of this case coincides with the $ d $-dimensional Euclidean norm $ |\cdot| $, and hence, the subdifferential $ \partial \gamma_0 $ coincides with the set valued function $ \bigl[ \frac{(\cdot)}{|\,\cdot\,|} \bigr] : \R^d \longrightarrow 2^{\mathbb{R}^d} $, defined in \eqref{Sgn^d}. 

    In the meantime, when $ \varepsilon > 0 $, the convex function $ \gamma_\varepsilon  $ belongs to $ C^\infty $-class, and the subdifferential $ \partial \gamma_\varepsilon $ is identified with the (single-valued) usual gradient:
    \begin{equation*}
        \nabla \gamma_\varepsilon : y \in \mathbb{R}^d \mapsto \nabla \gamma_\varepsilon(y) = \frac{y}{\sqrt{\varepsilon^2 +|y|^2}} \in \R^d.
    \end{equation*}
\end{ex}
\medskip

Next, we mention the notions of functional-convergences. 
 
\begin{defn}[Mosco--convergence: cf. \cite{MR0298508}]\label{Def.Mosco}
    Let $ H $ be an abstract Hilbert space. Let $ \Psi : H \longrightarrow (-\infty, \infty] $ be a proper, l.s.c., and convex function, and let $ \{ \Psi_n \}_{n =1}^\infty $ be a sequence of proper, l.s.c., and convex functions $ \Psi_n : H \longrightarrow (-\infty, \infty] $, $ n = 1, 2, 3, \dots $.  Then, it is said that $ \Psi_n \to \Psi $ on $ H $, in the sense of Mosco, as $ n \to \infty $, iff. the following two conditions are fulfilled:
\begin{description}
    \item[(\hypertarget{M_lb}{M1}) The condition of lower-bound:]$ \ds \varliminf_{n \to \infty} \Psi_n(\check{w}_n) \geq \Psi(\check{w}) $, if $ \check{w} \in H $, $ \{ \check{w}_n  \}_{n =1}^\infty \subset H $, and $ \check{w}_n \to \check{w} $ weakly in $ H $, as $ n \to \infty $; 
    \item[(\hypertarget{M_opt}{M2}) The condition of optimality:]for any $ \hat{w} \in D(\Psi) $, there exists a sequence \linebreak $ \{ \hat{w}_n \}_{n =1}^\infty  \subset H $ such that $ \hat{w}_n \to \hat{w} $ in $ H $ and $ \Psi_n(\hat{w}_n) \to \Psi(\hat{w}) $, as $ n \to \infty $.
\end{description}
\end{defn}
\begin{defn}[$ \Gamma $-convergence: cf. \cite{MR1201152}]\label{Def.Gamma}
    Let $ H $ be an abstract Hilbert space. Let $ \Psi : H \longrightarrow (-\infty, \infty] $ be a proper functional, and let $ \{ \Psi_n \}_{n = 1}^\infty $ be a sequence of proper functionals $ \Psi_n : H \longrightarrow (-\infty, \infty] $, $ n = 1, 2, 3, \dots $.  Then, it is said that $ \Psi_n \to \Psi $ on $ H $, in the sense of $ \Gamma $-convergence, as $ n \to \infty $, iff. the following two conditions are fulfilled:
\begin{description}
    \item[(\hypertarget{Gm_lb}{\boldmath$\Gamma$1}) The condition of lower-bound:]$ \ds \varliminf_{n \to \infty} \Psi_n(\check{w}_n) \geq \Psi(\check{w}) $, if $ \check{w} \in H $, $ \{ \check{w}_n  \}_{n = 1}^\infty \subset H $, and $ \check{w}_n \to \check{w} $ (strongly) in $ H $, as $ n \to \infty $; 
    \item[(\hypertarget{Gm_opt}{\boldmath$\Gamma$2}) The condition of optimality:]for any $ \hat{w} \in D(\Psi) $, there exists a sequence \linebreak $ \{ \hat{w}_n \}_{n = 1}^\infty  \subset H $ such that $ \hat{w}_n \to \hat{w} $ in $ H $ and $ \Psi_n(\hat{w}_n) \to \Psi(\hat{w}) $, as $ n \to \infty $.
\end{description}
\end{defn}

\begin{rem}\label{Rem.M-Gconvs}
Note that if the functionals are convex, then Mosco convergence implies $ \Gamma $-convergence, i.e., the $ \Gamma $-convergence of convex functions can be regarded as a weak version of Mosco convergence. Additionally, 
    in the Mosco-convergence of convex functions, we can see the following:
\begin{description}
    \item[(\hypertarget{Fact1}{Fact\,1})](cf. \cite[Theorem 3.66]{MR0773850} and \cite[Chapter 2]{Kenmochi81}) Let us assume that
    \begin{equation*}
        \Psi_n \to \Psi \mbox{ on $ H $, in the sense of  $\Gamma$-convergence, as $ n \to \infty $,}
    \vspace{-1ex}
\end{equation*}
and
\begin{equation*}
\left\{ ~ \parbox{10cm}{
$ [w, w^*] \in H \times H $, ~ $ [w_n, w_n^*] \in \partial \Psi_n $ in $ H \times H $, $ n \in \N $,
\\[1ex]
$ w_n \to w $ in $ H $ and $ w_n^* \to w^* $ weakly in $ H $, as $ n \to \infty $.
} \right.
\end{equation*}
Then, it holds that:
\begin{equation*}
[w, w^*] \in \partial \Psi \mbox{ in $ H \times H $, and } \Psi_n(w_n) \to \Psi(w) \mbox{, as $ n \to \infty $.}
\end{equation*}
    \item[(\hypertarget{Fact2}{Fact\,2})](cf. \cite[Lemma 4.1]{MR3661429} and \cite[Appendix]{MR2096945}) Let $ d \in \mathbb{N} $ denote the dimension constant, and let $  S \subset \R^d $ be a bounded open set. 
If: 
\begin{align*}
 \Psi_n \to \Psi \mbox{ on $ H $, in the sense of  Mosco, as $ n \to \infty $,}
\end{align*}
a sequence $ \{ \widehat{\Psi}_n^S \}_{n = 1}^\infty$ of proper, l.s.c., and convex functions on $ L^2(S; H) $, defined as:
        \begin{equation*}
            w \in L^2(S; H) \mapsto \widehat{\Psi}_n^S(w) := \left\{ \begin{array}{ll}
                    \multicolumn{2}{l}{\ds \int_S \Psi_n(w(t)) \, dt,}
                    \\[1ex]
                    & \mbox{ if $ \Psi_n(w) \in L^1(S) $,}
                    \\[2.5ex]
                    \infty, & \mbox{ otherwise,}
                \end{array} \right. \mbox{for $ n = 1, 2, 3, \dots $;}
        \end{equation*}
        converges to a proper, l.s.c., and convex function $ \widehat{\Psi}^S $ on $ L^2(S; H) $, defined as:
        \begin{equation*}
            z \in L^2(S; H) \mapsto \widehat{\Psi}^S(z) := \left\{ \begin{array}{ll}
                    \multicolumn{2}{l}{\ds \int_S \Psi(z(t)) \, dt, \mbox{ if $ \Psi(z) \in L^1(S) $,}}
                    \\[2ex]
                    \infty, & \mbox{ otherwise;}
                \end{array} \right. 
        \end{equation*}
        on $ L^2(S; H) $, in the sense of Mosco, as $ n \to \infty $. 
\end{description}
\end{rem}

\noindent
\underline{\textbf{\textit{Notations in BV-theory. (cf. \cite{MR1857292, MR2192832, MR3409135, MR0775682})}}}
Let $ d \in \N $ be a fixed constant of dimension. We denote by $\mathcal{L}^d$ the $d$-dimensional Lebesgue measure. The measure theoretical phrases, such as ``a.e.'', ``$dt$'', ``$dx$'', and so on, are all with respect to the Lebesgue measure in each corresponding dimension.

For any open set $ U \subset \mathbb{R}^d $, we denote by $ \mathcal{M}(U) $ (resp. $ \mathcal{M}_{\rm loc}(U) $) the space of all finite Radon measures (resp. the space of all Radon measures) on $ U $. In general, the space $ \mathcal{M}(U) $ (resp. $ \mathcal{M}_{\rm loc}(U) $) is known as the dual of the Banach space $ C_0(U) $ (resp. dual of the locally convex space $ C_{\rm c}(U) $), for any open set $ U \subset \mathbb{R}^d $. 

A function $ v \in L^1(U) $ (resp. $ v \in L_{\rm loc}^1(U) $)  is called a function of bounded variation, or a BV-function, (resp. a function of locally bounded variation or a BV$\empty_{\rm loc}$-function) on $ U $, iff. its distributional differential $ D v $ is a finite Radon measure on $ U $ (resp. a Radon measure on $ U $), namely $ D v \in \mathcal{M}(U) $ (resp. $ D v \in \mathcal{M}_{\rm loc}(U) $).
We denote by $ BV(U) $ (resp. $ BV_{\rm loc}(U) $) the space of all BV-functions (resp. all BV$\empty_{\rm loc}$-functions) on $ U $. For any $ v \in BV(U) $, the Radon measure $ D v $ is called the variation measure of $ v $, and its  total variation $ |Dv| $ is called the total variation measure of $ v $. Additionally, the value $|Dv|(U)$, for any $v \in BV(U)$, can be calculated as follows:
\begin{equation*}
|Dv|(U) = \sup \left\{ \begin{array}{l|l}
    \ds \int_{U} v \ {\rm div} \,\bm{\varpi} \, dx & \bm{\varpi} \in [C_\mathrm{c}^{1}(U)]^d \ \ \mbox{and}\ \ |\bm{\varpi}| \le 1\ \mbox{on}\ U
\end{array}
\right\}.
\end{equation*}
The space $BV(U)$ is a Banach space, endowed with the following norm:
\begin{equation*}
|v|_{BV(U)} := |v|_{L^{1}(U)} + |D v|(U),\ \ \mbox{for any}\ v\in BV(U).
\end{equation*}
We say that a sequence $\{ v_{n} \}_{n = 1}^\infty \subset BV(U)$ strictly converges in $BV(U)$ to $v \in BV(U)$ iff. $v_{n} \to v$ in $L^{1}(U)$ and $|Dv_{n}|(U) \to |Dv|(U)$ as $ n \to \infty $. In particular, if the boundary $\partial U$ is Lipschitz, then the space $BV(U)$ is continuously embedded into $L^{d/(d-1)}(U)$ and compactly embedded into $L^{q}(U)$ for any $1 \le q < d/(d-1)$ (cf. \cite[Corollary 3.49]{MR1857292} or \cite[Theorem 10.1.3-10.1.4]{MR2192832}). Additionally, if $1 \le r < \infty$, then the space $C^{\infty}(\overline{U})$ is dense in $BV(U) \cap L^{r}(U)$ for the \emph{intermediate convergence} (cf. \cite[Definition 10.1.3. and Theorem 10.1.2]{MR2192832}), i.e. for any $v \in BV(U) \cap L^{r}(U)$, there exists a sequence $\{v_{n} \}_{n = 1}^\infty \subset C^{\infty}(\overline{U})$ such that $v_{n} \to v$ in $L^{r}(U)$ and $\int_{U}|\nabla v_{n}|dx \to |Dv|(U)$ as $n \to \infty$.
\bigskip

\noindent
\underline{\textbf{\textit{Specific notations of this paper.}}} 
As is mentioned in the previous section, let $ (0, T) \subset \R$ be a bounded time-interval with a finite constant $ T > 0 $, and let $N \in \mathbb{N}$ be a constant of spatial dimension. Let $\Omega \subset \R^N $ be a fixed spatial bounded domain with a smooth boundary $ \Gamma := \partial\Omega $. We denote by $n_\Gamma$ the unit outward normal vector on $\Gamma$. Besides, we set $ Q := (0, T) \times \Omega $ and $ \Sigma := (0, T) \times \Gamma $. Especially, we denote by $ \partial_t $, $ \nabla $, and $ \mathrm{div} $ the distributional time-derivative, the distributional gradient, and distributional divergence, respectively. 
\medskip

On this basis, we define  
\begin{align*} 
    & \begin{cases}
        X := L^2(\Omega), ~~ Y := H^1(\Omega), 
        \\[1ex]
        W_{\nu} := \left\{z \in BV(\Omega) \cap X : \nu z \in Y \right\} \mbox{ for any } \nu \geq 0.
    \end{cases}
\end{align*} 

For any $ 0 \leq \beta \in Y \cap L^\infty(\Omega) $ and any $ z \in X $, we call the value $ {\rm Var}_\beta(z) \in [0, \infty] $, defined as,
\begin{equation*}
{\rm Var}_\beta(v) := \sup \left\{ \begin{array}{l|l}
    \displaystyle \int_\Omega v \, {\rm div} \, \bm{\varpi}\, dx & \parbox{6cm}{$ \bm{\varpi} \in [L^\infty(\Omega)]^N $ with a compact support, and $ | \bm{\varpi} | \leq \beta $ a.e.\ in $ \Omega $}
\end{array} \right\} \in [0, \infty],
\end{equation*}
``the total variation of $ v $ weighted by $ \beta $'', or the ``weighted total variation'' in short. 

\begin{rem}\label{Rem.Note05}
Referring to the general theories (e.g., \cite{MR1259102,MR1857292,MR1736243}), we can confirm the following facts associated with the weighted total variations. 
\begin{description}
	\item[{({Fact\,3})}](cf.\ \cite[Theorem 5]{MR1736243}) For any $ 0 \leq \beta \in Y \cap L^\infty(\Omega) $, the functional $ z \in X \mapsto {\rm Var}_\beta(z) \in [0, \infty] $ is a proper, l.s.c.,\ and convex function that coincides with the lower semi-continuous envelope of
\begin{equation*}
z \in W^{1, 1}(\Omega) \cap X \mapsto \int_\Omega \beta |D z| \, dx \in [0, \infty).
\end{equation*}
\item[{({Fact\,4})}](cf.\ \cite[Theorem 4.3]{MR1259102} and \cite[Proposition 5.48]{MR1857292}) If $ 0 \leq \beta \in Y \cap L^\infty(\Omega) $ and $ z \in BV(\Omega) \cap L^2(\Omega) $, then there exists a Radon measure $ |Dz|_\beta \in \mathcal{M}(\Omega) $ such that
\begin{equation*}
|Dz|_\beta(\Omega) = \int_\Omega  d|Dz|_\beta = {\rm Var}_\beta(z),
\end{equation*}
and 
\begin{equation}\label{|Dz|_beta(A)}
\begin{array}{c}
\left\{ ~ \parbox{12.25cm}{
$ |Dz|_\beta(A) \leq |\beta|_{L^\infty(\Omega)} |Dz|(A) $,
\\[1ex]
$ \displaystyle
|Dz|_\beta(A) = \inf \left\{ \begin{array}{l|l}
		\displaystyle \liminf_{n \to \infty} \int_A \beta |D \tilde{z}_n| \, dx & \parbox{5cm}{$ \{ \tilde{z}_n \}_{n = 1}^\infty \subset W^{1, 1}(A) \cap L^2(A) $ such that \ $ \tilde{z}_n \to z $ in $ L^2(A) $  as $ n \to \infty $}
\end{array} \right\} $,
} \right.
\end{array}
\end{equation}
for any open set $ A \subset \Omega $.
\item[{({Fact\,5})}]If $ \beta \in Y $, $ \log \beta \in L^\infty(\Omega) $, and $ z \in BV(\Omega) \cap X $, then it follows that
\begin{equation}\label{|Dz|_BV}
\left\{ ~ {\parbox{10cm}{
$ |D{z}|_\beta(A) \geq c_\beta |Dz|(A) $ \ for any open set $ A \subset \Omega $,
\\[2ex]
$ D({\rm Var}_\beta) = BV(\Omega) \cap X $, \ and
\\[1ex]
$ {\rm Var}_\beta(z) = \sup \left\{ \begin{array}{l|l}
    \displaystyle \int_\Omega z \, {\rm div} \, (\beta \bm{\varpi}) \, dx & \parbox{4cm}{$ \bm{\varpi} \in [L^\infty(\Omega)]^N $ with a  compact support, and $ | \bm{\varpi} | \leq 1 $ a.e.\ in $ \Omega $}
\end{array} \right\}, $}} \right.
\end{equation}
\end{description}
where $C_{\beta} := \mathrm{essinf}_{x \in \Omega} \beta(x) > 0$.
Moreover, the following properties can be inferred from \eqref{|Dz|_beta(A)}--\eqref{|Dz|_BV}:
\begin{itemize}
\item[$\bullet$] $ |Dz|_c = c|Dz| $ in $ \mathcal{M}(\Omega) $ for any constant $c \geq 0 $ and $ z \in BV(\Omega) \cap X $;
\item[$\bullet$] $ |Dz|_\beta = \beta |D z| \mathscr{L}^N $ in $ \mathcal{M}(\Omega) $, \ if $ 0 \leq \beta \in H^1(\Omega) \cap L^\infty(\Omega) $ and $ z \in W^{1, 1}(\Omega) \cap X $.
\end{itemize}
\end{rem}

\begin{defn}[Generalized weighted total variation; cf.\ {\cite[Section 2]{MR3268865}}]\label{Note06}
    \ For any (possibly sign-changing) $ \beta \in Y \cap L^\infty(\Omega) $ and any $ z \in BV(\Omega) \cap X $, we define a real-valued Radon measure $ [\beta |Dz|] \in \mathcal{M}(\Omega) $, as follows:
\begin{equation*}
[\beta |Dz|](B) := |Dz|_{[\beta]^+}(B) - |Dz|_{[\beta]^-}(B) \mbox{ \ for any Borel set $ B \subset \Omega $.}
\end{equation*}
Note that $ [\beta|D z|](\Omega) $ can be configured as a {\em generalized total variation} of $ z \in BV(\Omega) \cap X $ by the possibly sign-changing weight $ \beta \in Y \cap L^\infty(\Omega) $. 
\end{defn}

\begin{rem}\label{Rem.Note06}
With regard to the generalized weighted total variations, the following facts are verified in \cite[Section 2]{MR3268865}.
\begin{description}
\item[{({Fact\,6})}]
Let $ \beta \in Y \cap L^\infty(\Omega) $ and $ z \in BV(\Omega) \cap X $ be arbitrary fixed functions, and let $ \{ z_n \}_{n = 1}^\infty \subset C^\infty(\overline{\Omega}) $ be a sequence such that
\begin{equation*}
z_n \to z \mbox{ in $ X $ \ and \ strictly \ in } BV(\Omega) \mbox{ as $ n \to \infty $.}
\end{equation*}
Then
\begin{equation*}
\int_\Omega \beta |D z_n| \, dx \to \int_\Omega d[\beta |Dz|] \mbox{ \ as $ n \to \infty $.}
\end{equation*}
\item[{({Fact\,7})}]For any $ z \in BV(\Omega) \cap X $, the mapping
\begin{equation*}
\ds \beta \in H^1(\Omega) \cap L^\infty(\Omega) \mapsto \int_\Omega d[\beta |Dz|] \in \R
\end{equation*}
is a linear functional, and moreover, if $ \varphi \in H^1(\Omega) \cap C(\overline{\Omega}) $ and $ \beta \in Y \cap L^\infty(\Omega) $, then
\begin{equation*}
\displaystyle \int_\Omega d[\varphi \beta |Dz|] = \int_\Omega \varphi \, d[\beta |Dz|].
\end{equation*}
\end{description}
\end{rem}

Finally, we recall the previous works on time-dependent total variation, obtained in \cite{MR3268865,MR3462536}.
\begin{prop}\textbf{(cf. \cite[Lemma 5]{MR3268865} and \cite[Remark 2]{MR3462536})}\label{PropKen01}
    Let us fix an open interval $ I \subset (0, T) $, and let us take any function $ \zeta \in L^2(I; X) $ satisfying $ |D\zeta(\cdot)|(\Omega) \in L^1(I) $. Then, There exists $ \{ \omega_n \}_{n = 1}^\infty \subset C^\infty(\overline{I \times \Omega}) $, such that:
        \begin{equation*}
            \omega_n \to \zeta \mbox{ in $ L^2(I; X) $, and } \int_I \left| \int_\Omega |\nabla \omega_n| \, dx -|D \zeta(t)|(\Omega) \right| \, dt \to 0, \mbox{ as $ n \to \infty $.}
        \end{equation*}
\end{prop}
\begin{prop}\label{rem-2}\emph{(cf. \cite[Lemmas 6 and 7]{MR3268865})}
    Let us fix an open interval $ I \subset (0, T) $, and let us assume 
    \begin{equation*}
\left\{ ~ \parbox{12cm}{
		$ \beta \in C(\overline{I}; X) \cap L^\infty(I; Y) $, \ $ \log \beta \in L^\infty(I \times \Omega) $,
	\\[1ex]
        $ \{ \beta_n \} \subset  L^2(I; Y) $, and $ \beta_n \geq 0 $ a.e. in $ I \times \Omega $,  for $ n = 1,2, 3, \dots $,
	\\[1ex]
	$ \beta_n(t) \to \beta(t) $ in $ X $ and weakly in $ Y $ as $ n \to \infty $, for a.e. $ t \in I $,
} \right.
\end{equation*}
and
\begin{equation*}
\left\{ ~ \hspace{-5ex} \parbox{10cm}{
\vspace{-2ex}
\begin{itemize}
    \item[]$ \zeta \in C(\overline{I}; X) $, ~
        $ \{ \zeta_n \}_{n = 1}^\infty \subset L^2(I; Y) $, 
\vspace{-1ex}
\item[]$ \zeta_n(t) \to \zeta(t) $ in $ X $ as $ n \to \infty $, \ a.e.\ $ t \in I $.
\vspace{-2ex}
\end{itemize}
} \right.
\end{equation*}
   Then, the following items hold.
\begin{description}
    \item[\textmd{(I)}]The function:
        $ t \in I \mapsto [\beta(t)|D \zeta(t)|](\Omega) \in [0, \infty] $ 
        is lower semi-continuous (possibly not proper), and
\begin{equation*}
\liminf_{n \to \infty}\int_I \int_\Omega \beta_n(t) |\nabla \zeta_n(t)| \, dx dt \geq \int_I \int_\Omega d[\beta(t) |D \zeta(t)|] \, dt. 
\end{equation*}
\item[\textmd{(II)}]If:
\begin{equation*}
\int_I \int_\Omega d [\beta_n(t) |D \zeta_n(t)|] \, dt \to \int_I \int_\Omega d[\beta(t) |D \zeta(t)|] \, dt \mbox{ as $ n \to \infty $} 
\end{equation*}
and
\begin{equation*}
\left\{ \hspace{-4ex} \parbox{15cm}{
\vspace{-1ex}
\begin{itemize}
\item[]$ \varrho \in L^\infty(I; H^1(\Omega)) \cap L^\infty(I \times \Omega) $, $ \{ \varrho_n \}_{n = 1}^\infty \subset L^\infty(I; H^1(\Omega)) \cap L^\infty(I \times \Omega) $,
\vspace{0ex}
\item[]$ \{ \varrho_n \}_{n = 1}^\infty $ is a bounded sequence in $ L^\infty(I \times \Omega) $,
\vspace{0ex}
\item[]$ \varrho_n(t) \to \varrho(t) $ in $ L^2(\Omega) $ and weakly in $ H^1(\Omega) $ as $ n \to \infty $, a.e.\ $ t \in I $,
\vspace{-1ex}
\end{itemize}
} \right.
\end{equation*}
then
\begin{equation*}
\int_I \int_\Omega \varrho_n(t) |\nabla \zeta_n(t)| \, dx dt \to \int_I \int_\Omega d [\varrho(t) |D \zeta(t)|] \mbox{ as $ n \to \infty $.}
\end{equation*}
\end{description}
\end{prop}

\section{Main Theorem}
We begin by setting up the assumptions required for our Main Theorem. 
We discuss the Main Theorem under the following assumptions:  
\begin{description}
\item[\textmd{(\hypertarget{A1l}{A1})}]
    $\kappa > 0$, $M_0 > 0$, and $\nu_0 \geq 0$ are fixed constants.
\item[\textmd{(\hypertarget{A2l}{A2})}]
    $  [u, v]  \in [L^\infty(Q)]^2$ is a fixed pair of functions.
\item[\textmd{(\hypertarget{A3l}{A3})}]
    $\alpha_{0} \in W_\mathrm{loc}^{1, \infty}(\R)$ is a fixed function with the first derivative $ \alpha_0' = \frac{d \alpha}{d \eta} \in L_\mathrm{loc}^\infty(\R) $, $\alpha \in C^2(\R)$ is a fixed  with the first derivative $ \alpha' = \frac{d \alpha}{d \eta} \in C^1(\R) $ and second derivative $ \alpha'' = \frac{d^2 \alpha}{d\eta^2} \in C(\R) $. Furthermore, $ \alpha_0 $ and $ \alpha $ fulfill the following:
        \begin{itemize}
            \item $\alpha'(0) = 0$, $ \alpha'' \geq 0 $ on $ \R $;
            \item $ \inf \alpha(\R) \cup \alpha_0(\R) \geq \delta_* $ for some constant $ \delta_* \in (0, 1) $.
        \end{itemize}
    \item[\textmd{(\hypertarget{A4l}{A4})}]
Let $g : \mathbb{R} \longrightarrow  \mathbb{R}$ be a $C^{2}$-function, which is Lipschitz continuous on $\mathbb{R}$. Additionally, $g$ has a nonnegative primitive $ 0 \leq G \in C^{3}(\mathbb{R})$; that is, the derivative $ G'= \frac{dG}{d\eta} $ coincides with $ g $ on $ \R $. Moreover, $g$ satisfies the following:
        \begin{equation*}
            \liminf_{\xi \downarrow -\infty}g(\xi) = -\infty \mbox{ and } \limsup_{\xi \uparrow \infty}g(\xi) = \infty.
        \end{equation*}
\end{description}

Using this, we define the time-periodic solution to (S) as follows: 

\begin{rem}\label{Rem.bdd}
    From (\hyperlink{A2l}{A2}) and (\hyperlink{A4l}{A4}), we immediately determine a constant $ R_0 > 0 $ such that
    \begin{align}\label{ken00}
        & \begin{cases}
            |u|_{L^\infty(Q)} \vee |v|_{L^\infty(Q)} \leq R_0,  
            \\[1ex]
            g(-R_0) \leq -|u|_{L^\infty(Q)} \leq |u|_{L^\infty(Q)} \leq g(R_0),
            \\[1ex]
            M_0 R_0 \geq |v|_{L^\infty(Q)}. 
        \end{cases}
\end{align}
\end{rem}

\begin{defn}[T-periodic solution]\label{Def.peri} 
    A function $[\eta, \theta] \in L^2(0, T; [X]^2)$ is called a time-periodic solution to (S) on $ (0, T) $, or \textit{T-periodic solution} in short, iff
        \begin{equation}\label{S1}
            \begin{cases}
                [\eta, \theta ] \in  W^{1, 2}(0, T; [X]^2) \cap [L^\infty(Q)]^2, 
                \\[1ex]
                [\eta, \nu_0 \theta] \in [L^\infty(0, T; Y)]^2, ~~ |D\theta(\cdot)|(\Omega) \in L^\infty(0, T),
                \\[1ex]
                 [\eta(0), \theta(0)] = [\eta(T), \theta(T)] \mbox{ in } [X]^2,
       \end{cases}
        \end{equation}
        \begin{equation}\label{S2}
                \begin{array}{c}
                    \ds\bigl( \partial_{t}\eta(t), \varphi \bigr)_{X} +\kappa^2 \bigl( \nabla \eta(t), \nabla \varphi \bigr)_{[X]^N} +\bigl( g(\eta(t)), \varphi \bigr)_{X} 
                    \\[1.5ex]
                    \ds + \int_\Omega d \bigl[\varphi \alpha(\eta(t))|D\theta(t)| \bigr]\,  = \bigl(  u(t), \varphi \bigr)_{X} 
                    \\[2ex]
                    \mbox{ for any $ \varphi \in Y \cap L^\infty(\Omega) $, a.e. $ t \in (0, T) $;}
                \end{array}
        \end{equation}
and
        \begin{equation}\label{S3}
                \begin{array}{c}
                    \ds\bigl( \alpha_{0}(\eta(t))\partial_{t}\theta(t), \theta(t)-\psi \bigr)_{X}  +M_0 \bigl( \theta(t), \theta(t) -\psi \bigr)_{X}
                    \\[1ex]
                    \ds + \bigl( \nabla (\nu_0 \theta)(t), \nabla \bigl( \nu_0(\theta(t)-\psi) \bigr) \bigr)_{[X]^N} +\int_{\Omega} d \bigl[ \alpha(\eta(t)|D\theta(t)| \bigr] 
                    \\[1ex]
                    \ds \leq \int_{\Omega} d \bigl[ \alpha(\eta(t))|\nabla\psi| \bigr]
                     +\bigl(  v(t), \theta(t)-\psi \bigr)_{X} 
                     ~\mbox{ for any } \psi \in W_{\nu_0}, \mbox{ a.e. } t \in (0, T). 
                \end{array}
        \end{equation}
\end{defn}
\bigskip

Now, our goal in this paper is to prove the Main Theorem.
\paragraph{Main Theorem.}{\em
We assume (\hyperlink{A1l}{A1})--(\hyperlink{A4l}{A4}). 
Then, the system (S) admits at least one T-periodic solution $[\eta, \theta] \in L^2(0, T; [X]^2)$.
}

\begin{rem}\label{Rem.1stEq}
    We note that the notion of weighted total variation is necessary for the rigorous definition of the $ T $-periodic solution. This notion allows the definition formula \eqref{S2} to be nonstandard as a parabolic variational inequality. However, if we consider a sufficiently large exponent $ s > N/2 $, then (\hyperlink{Fact7l}{Fact\,7}) and the embedding $ H^{s}(\Omega) \subset C(\overline{\Omega}) $ enable us to reformulate \eqref{S2} as an evolution equation on a Hilbert space $ H^{s}(\Omega)^* $:
    \begin{align*}
        \partial_t \eta(t) & ~+\kappa^2 \bigl( F_Y^* \eta(t) -\eta(t) \bigr) +g(\eta(t)) +[\alpha'(\eta(t))|D \theta(t)|] 
        \\
        & = u(t) \mbox{ in $ H^{s}(\Omega)^* $, a.e. $ t \in (0, T) $,} 
    \end{align*}
    where $ F_Y^* : Y \longrightarrow Y^* $ $ (\subset H^{s}(\Omega)^*) $ is the duality mapping between $ Y $ and $ Y^* $; that is,
    \begin{equation*}
        \langle F_Y^* \tilde{\eta}, \varphi \rangle_Y := (\nabla \tilde{\eta}, \nabla \varphi)_{[X]^N} +(\tilde{\eta}, \varphi)_X, ~\mbox{for all $ [\tilde{\eta}, \varphi] \in [Y]^2 $.}
    \end{equation*}
\end{rem}

\begin{rem}\label{Rem.2ndEq}
    We consider any $ \nu \geq 0 $ and define a functional $ \Phi_\nu : [X]^2 \longrightarrow [0, \infty] $ as follows:
\begin{align}\label{phi.nu}
    \displaystyle \displaystyle [\tilde{\eta}, \tilde{\theta}]  \in [X]^2 \mapsto \Phi_\nu (\tilde{\eta}, \tilde{\theta}) := \left\{
\begin{array}{l}
    \displaystyle 
    \int_\Omega d \bigl[ \alpha(\tilde{\eta})|D \tilde{\theta}| \bigr] +\frac{1}{2}\int_{\Omega} |\nabla (\nu\tilde{\theta})|^2 \, dx 
\\[2.0ex]
    \qquad \ds  +\frac{M_0}{2} \int_\Omega |\tilde{\theta}|^2 \, dx,\mbox{ if } \tilde{\theta} \in W_\nu, \\[4ex]
\infty, \mbox{ otherwise.}
\end{array}
\right.
\end{align}
    Then, because of Proposition \ref{rem-2}, we observe that $ \Phi_\nu(\tilde{\eta}, \cdot) : X \longrightarrow [0, \infty] $ is a proper, l.s.c., and convex function for each fixed $ \tilde{\eta} \in Y \cap L^\infty(\Omega) $. Given this, we can reformulate the variational inequality \eqref{S3} as follows:
    \begin{equation*}
        \alpha_0(\eta(t)) \partial_t \theta(t) +\partial \Phi_\nu(\eta(t), \theta(t)) \ni v(t) \mbox{ in $ X $, a.e. $ t \in (0, T) $}
    \end{equation*}
    using the subdifferentials $ \partial \Phi_\nu(\eta(t), \cdot) \subset [X]^2 $ of the time-dependent convex functions $ \Phi_\nu(\eta(t), \cdot) : X \longrightarrow [0, \infty] $ for every $ t \in (0, T) $. 

    Moreover, it should be noted that the rigorous definition of free energy$ \mathcal{F}_\nu $, in \eqref{freeEnergyOrg}, is given by
    \begin{align}
        \mathcal{F}_\nu : \, & [\tilde{\eta}, \tilde{\theta}] \in D(\mathcal{F}_\nu) := D(\Phi_\nu) \cap [Y \times W_\nu] \subset [X]^2
        \nonumber
        \\
        & \mapsto \mathcal{F}_\nu(\tilde{\eta}, \tilde{\theta}) := \frac{\kappa^2}{2} \int_\Omega |\nabla \tilde{\eta}|^2 \, dx +\int_\Omega G(\eta) \, dx +\Phi_\nu(\tilde{\eta}, \tilde{\theta}) \in [0, \infty).
        \label{freeEnergyReal}
    \end{align}
\end{rem}

\section{Time-discrete approximation} 
    In this section, we consider a system of time-discretization schemes as the approximating problem of our system (S). We consider the time-discrete approximation under the following assumptions. 

\begin{description}
    \item[\textmd{(\hypertarget{A5l}{A5})}]The time interval $ [0, T] $ is discretized by a finite number of division  points:
        \begin{align*}
            t_i := i \tau, ~& i = 0, 1, 2, 3, \dots, m, ~ \mbox{ with division size $ m \in \N $,}
            \\
            & \mbox{and time-step size }~ \tau := \frac{T}{m}. 
        \end{align*}
        Additionally, division size $ m $ is sufficiently large to satisfy the following:   
    \begin{align*}
        m > 4 T \bigl( |g'|_{L^\infty(\mathbb{R})} + 1 \bigr); \mbox{hence, } 0 < \tau < \tau^* := \frac{1}{4 \bigl( |g'|_{L^\infty(\mathbb{R})} + 1 \bigr)}.
    \end{align*}
\item[\textmd{(\hypertarget{A6l}{A6})}]
    The functional $ \Phi_0 $, given in \eqref{phi.nu}, is relaxed by a sequence of functionals $ \{ \Phi_{\nu, \varepsilon} \}_{\nu, \varepsilon \in (0, 1)} $, which is defined as follows:
\begin{align}\label{phi.nu-eps}
    & \displaystyle \displaystyle [\tilde{\eta}, \tilde{\theta}]  \in [X]^2 \mapsto \Phi_{\nu, \varepsilon} (\tilde{\eta}, \tilde{\theta}) := \left\{
\begin{array}{l}
\displaystyle \int_\Omega \alpha(\tilde{\eta})\gamma_\varepsilon(\nabla \tilde{\theta}) \, dx +\frac{\nu^2}{2}\int_{\Omega} |\nabla \tilde{\theta}|^2 \, dx, 
\\[2.0ex]
    \, \, \quad \displaystyle +\frac{M_0}{2} \int_\Omega |\tilde{\theta}|^2 \, dx, \mbox{ if } \tilde{\theta} \in Y, \\[1.5ex]
\infty, \mbox{ otherwise,}
\end{array}
\right.
        \\[2ex]
        & \mbox{for every $ \nu > 0 $ and $ \varepsilon \in (0, 1) $.}
        \nonumber
\end{align}
Note that $\gamma_\varepsilon$ is the continuous and convex function as in \eqref{gamma_eps}.
\item[\textmd{(\hypertarget{A7l}{A7})}]
    Let $[\{ u_i^{m}\}_{i=1}^m, \{v_i^{m} \}_{i=1}^m] \in  [X]^{m \times 2} $ be the time-discretization approximation of the forcing pair $ [u, v] \in [L^\infty(Q)]^2 $, which is defined as follows: 
\begin{align*}
    \displaystyle [u_i^{m}, v_i^{m}] :=  \left[ \frac{1}{\tau} \int_{t_{i -1}}^{t_i}u(t)\, dt,\, \frac{1}{\tau} \int_{t_{i-1}}^{t_i}v(t)\, dt  \right] \mbox{ in } [X]^2 \mbox{ for any } i = 1, 2, 3, \ldots, m.
\end{align*}
\end{description}

\begin{rem}\label{uv.P}
    Clearly, from  (\hyperlink{A2l}{A2}),  (\hyperlink{A7l}{A7}), and Remark \ref{Rem.bdd},
    \begin{align*}
    [u_i^{m}, v_i^{m}] \in Z(R_0) :=  ~&\left\{ \begin{array}{l|l}
        [\tilde{u}, \tilde{v}] \in [L^\infty(\Omega)]^2 & \hspace{-3ex}\parbox{5.6cm}{
            \vspace{-2ex}
            \begin{itemize}
                \item $ |\tilde{u}|_{L^\infty(\Omega)} \vee |\tilde{v}|_{L^\infty(\Omega)} \leq R_0 $, 
                    \vspace{-1ex}
                \item $ g(-R_0) \leq -|\tilde{u}|_{L^\infty(\Omega)} $, 
                    \vspace{-1ex}
                \item $ |\tilde{u}|_{L^\infty(\Omega)} \leq g(R_0) $, 
                    \vspace{-1ex}
                \item $ M_0 R_0 \geq |\tilde{v}|_{L^\infty(\Omega)} $
            \vspace{-2ex}
            \end{itemize}
        }
    \end{array} \right\}
    \\[1ex]
    & \mbox{for $ i = 1, 2, 3, \dots, m $,}
    \nonumber
\end{align*}
    where $ R_0 > 0 $ is the constant as in \eqref{ken00}. Additionally, we note that the sequence $ \{ \Phi_{\nu, \varepsilon} \}_{\varepsilon \in (0, 1)} $, as in (\hyperlink{A6l}{A6}), results in the following relaxation sequence $ \{ \mathcal{F}_{\nu, \varepsilon} \}_{\varepsilon \in (0, 1)} $ for free energy$ \mathcal{F}_\nu $, given in \eqref{freeEnergyReal}, which is defined as follows:  
    \begin{align}\label{free.en}
        [\tilde{\eta}, \tilde{\theta}]  \in [X]^2 \mapsto & \displaystyle \mathcal{F}_{\nu, \varepsilon} (\tilde{\eta}, \tilde{\theta}) := \, \frac{\kappa^2}{2} \int_\Omega |\nabla \tilde{\eta}|^2 \, dx  + \int_\Omega G(\tilde{\eta})\, dx 
        \\
        & +\Phi_{\nu, \varepsilon} (\tilde{\eta}, \tilde{\theta}) \mbox{ for $ \varepsilon \in (0, 1) $.}
        \nonumber
\end{align} 
\end{rem}

Now, for every $ \nu > 0 $, $ \varepsilon \in (0, 1) $, and $ \tau \in (0, \tau^*) $ as in (\hyperlink{A5l}{A5})--(\hyperlink{A7l}{A7}),  we let \hyperlink{(S)$_{\nu, \varepsilon}^{m}$}{(S)$_{\nu, \varepsilon}^{m}$} denote  the approximating problem of our system (S), which is described by the following system of time-discretization schemes:

~~ \hypertarget{(S)$_{\nu, \varepsilon}^{m}$}{(S)$_{\nu, \varepsilon}^{m}$}
\begin{align*}
    & \begin{cases}
        \displaystyle \frac{1}{\tau} \alpha_0(\eta_{i -1})(\theta_i - \theta_{i-1} ) - \mathrm{div} \left( \alpha(\eta_{i -1}) \nabla \gamma_\varepsilon(\nabla\theta_i) +\nu^2 \nabla \theta_i \right) +M_0\theta_i = v_i^{m} \mbox{ a.e. in } \Omega,
               \\[2ex]
               \displaystyle \left( \alpha(\eta_{i -1}) \nabla \gamma_\varepsilon(\nabla\theta_i) +\nu^2 \nabla \theta_i \right) \cdot n_\Gamma = 0 \mbox{ on }  \Gamma,
    \end{cases}
    \\[1ex]
    & \begin{cases}
        \displaystyle \frac{1}{\tau} (\eta_i - \eta_{i-1}) - \kappa^2 \mathit{\Delta} \eta_i +g(\eta_i) +\alpha'(\eta_i) \gamma_\varepsilon(\nabla \theta_{i}) =u_i^{m}  ~\mbox{ a.e. in } \Omega,
        \\[2ex]
        \nabla \eta_i \cdot n_\Gamma = 0  \mbox{ on } \Gamma 
    \end{cases}
    \\[1ex]
    & \mbox{for $ i = 1, 2, 3, \dots, m $, starting from the initial data $ [\eta_0, \theta_0] \in [Y]^2 $.}
 \end{align*}
\begin{defn}
    [Approximating solutions]
    \label{Def.time-dis}
A pair of sequences $ [\{\eta_i \}_{i=1}^m, \{ \theta_i\}_{i=1}^m ] \in [X]^{m \times 2}$ is called a solution to the approximating system \hyperlink{(S)$_{\nu, \varepsilon}^{m}$}{(S)$_{\nu, \varepsilon}^{m}$}, or the approximating solution in short, iff
\begin{subequations}\label{3tau}
    \begin{equation}\label{3tau-0}\noeqref{3tau-0}
        [\{\eta_i \}_{i=1}^m, \{ \theta_i\}_{i=1}^m ] \in  [Y]^{m \times 2} ,
    \end{equation}
    \vspace{-3ex}
\begin{align}\label{3tau-1}
    \displaystyle \frac{1}{\tau} (\eta_i - \eta_{i-1}& \displaystyle , \varphi)_X + \kappa^2( \nabla \eta_i, \nabla \varphi)_{[X]^N} + \bigl( g(\eta_i) +\alpha'(\eta_i) \gamma_\varepsilon(\nabla \theta_{i}), \varphi \bigr)_X 
    \\[0ex]
    & = (u_i^{m}, \varphi)_X ~\mbox{for any $\varphi \in Y $ and $ i = 1, 2, 3, \ldots, m  $,}
    \nonumber
\end{align}
    \vspace{-3ex}
\begin{align}\label{3tau-2}
    & \displaystyle  \frac{1}{\tau} (\alpha_{0}(\eta_{i -1})( \theta_i  - \theta_{i-1}), \psi)_X + \bigl( \alpha(\eta_{i -1})\nabla \gamma_\varepsilon (\nabla \theta_i) +\nu^2 \nabla\theta_i, \nabla \psi \bigr)_{[X]^N} 
    \\[0ex] 
    & + M_0(\theta_i, \psi)_X = (v_i^{m}, \psi)_X ~\mbox{ for any $\psi \in Y$ and  $ i = 1, 2, 3, \ldots, m$.}
    \nonumber
\end{align}
\end{subequations}
Additionally, an approximating solution $ [\{\eta_i \}_{i=1}^m, \{ \theta_i\}_{i=1}^m ] \in [X]^{m \times 2} $ is called a periodic solution, or approximating periodic solution, iff
    \begin{equation}\label{3tau-3}
        [\eta_m, \theta_m] = [\eta_0, \theta_0] \mbox{ in $ [X]^2 $.}
    \end{equation}
\end{defn}

\begin{rem}\label{Shorem01}
From assumption (\hyperlink{A6l}{A6}), we observe that \eqref{3tau-2} of Definition \ref{Def.time-dis} is equivalent to the following variational inequality: 
\begin{align*}
 \displaystyle  \left( \frac{1}{\tau} (\alpha_{0}(\eta_{i -1})( \theta_i  - \theta_{i-1})) - v_i^{m} , \theta_i - \psi \right)_X & \, + \Phi_{\nu, \varepsilon}(\eta_{i-1}, \theta_i) - \Phi_{\nu, \varepsilon}(\eta_{i-1}, \psi) \leq 0 \nonumber \\
 & \mbox{ for any } \psi \in Y \mbox{ and } i = 1, 2, 3, \ldots, m. 
\end{align*}
\end{rem}

Based on this, our goal in this section is to prove Theorem \ref{Thm.1}. 
\begin{thm}[Existence of an approximating periodic solution]\label{Thm.1} Under assumptions \linebreak (\hypertarget{A1l}{A1})--(\hyperlink{A7l}{A7}), the approximating system (S)$_{\nu, \varepsilon}^{m}$ admits at least one periodic solution \linebreak $ [\{\eta_i \}_{i=1}^m, \{ \theta_i\}_{i=1}^m ] \in  [X]^{m \times 2} $. 
\end{thm}

Before we prove Theorem \ref{Thm.1}, we prepare some lemmas. 
\begin{lem}\label{Lem3-02}
    Under the assumptions  and notation in  (\hyperlink{A1l}{A1})--(\hyperlink{A6l}{A6}) and Remark \ref{uv.P}, let $ [\bar{\eta}_0, \bar{\theta}_0] \in [Y]^2 $ and $ [\bar{u}, \bar{v}] \in [X]^2 $ be a fixed pairs of functions. 
Then, it holds that: 
 \begin{description}
     \item[\textmd{(\hypertarget{I}{I})}] A variational identity
            \begin{align}\label{ken01}
                \displaystyle \frac{1}{\tau} (\eta - \bar{\eta}_{0}, \varphi)_X + \kappa^2& \displaystyle ( \nabla \eta, \nabla \varphi)_{[X]^N} + \bigl( g(\eta) +\alpha'(\eta) \gamma_\varepsilon(\nabla \bar{\theta}_0), \varphi \bigr)_X 
                \\[0ex]
                & = (\bar{u}, \varphi)_X ~\mbox{for any $\varphi \in Y $}
                \nonumber
            \end{align}
            admits a unique solution $ \eta \in Y $. 
        \item[\textmd{(\hypertarget{II}{II})}] A  variational identity
        \begin{align}\label{ken02}
            \displaystyle  \frac{1}{\tau} (\alpha_{0}& \displaystyle (\bar{\eta}_0)( \theta  - \bar{\theta}_{0}), \psi)_X + \bigl( \alpha(\bar{\eta}_0)\nabla \gamma_\varepsilon (\nabla \theta) +\nu^2 \nabla\theta, \nabla \psi \bigr)_{[X]^N} 
            \\[0ex] 
            & + M_0(\theta, \psi)_X = (\bar{v}, \psi)_X ~\mbox{ for any $\psi \in Y$}
            \nonumber
        \end{align}
         admits a unique solution $ \theta \in Y $. 
     \item[\textmd{(\hypertarget{III}{III})}] In the variational identities \eqref{ken01} and \eqref{ken02}, if
         \begin{subequations}\label{ken03}
         \begin{equation}\label{ken03a}
             \begin{cases}
                 [\bar{u}, \bar{v}] \in Z(R_0), ~ 
                 [\bar{\eta}_0, \bar{\theta}_0] \in [L^\infty(\Omega)]^2, 
                 \\
                 \mbox{ and } |\bar{\eta}_0|_{L^\infty(\Omega)} \vee |\bar{\theta}_0|_{L^\infty(\Omega)} \leq R_0,
             \end{cases}
         \end{equation}
         then the corresponding solutions $ \eta $ and $ \theta $ fulfill the following: 
         \begin{equation}\label{ken03b}
            [\eta, \theta] \in [L^\infty(\Omega)]^2, \mbox{ and } |\eta|_{L^\infty(\Omega)} \vee |\theta|_{L^\infty(\Omega)} \leq R_0,
         \end{equation}
         \end{subequations}
         where $ R_0 > 0 $ is the constant as in \eqref{ken00} and $ Z(R_0) \subset [L^\infty(\Omega)]^2 $ is the class of functional pairs as in Remark \ref{uv.P}.  
\end{description}
\end{lem}
\paragraph{Proof of Lemma \ref{Lem3-02}.}{
    To prove item (I), we first recall assumptions (\hypertarget{A1l}{A1})--(\hyperlink{A6l}{A6}) to note that the variational identity \eqref{ken01} corresponds to the Euler--Lagrange equation for the proper and l.s.c. potential functional $ \Psi_0 : X \longrightarrow [0, \infty] $, which is defined as follows:
    \begin{equation*}
        \tilde{\eta} \in X \mapsto \Psi_0(\tilde{\eta}) := \begin{cases}
            \ds \frac{1}{2 \tau} |\tilde{\eta} -\bar{\eta}_0|_X^2 +\frac{\kappa^2}{2} \int_\Omega |\nabla \tilde{\eta}|^2 \, dx +\int_\Omega G(\tilde{\eta}) \, dx 
            \\[2ex]
            \qquad \ds +\int_\Omega \alpha(\tilde{\eta}) \gamma_\varepsilon(\nabla \bar{\theta}_0) \, dx -(\bar{u}, \tilde{\eta})_X, ~\mbox{ if $ \tilde{\eta} \in Y $,}
            \\[3ex]
            \infty, ~\mbox{ otherwise.}
        \end{cases}
    \end{equation*} 
    Particularly, from (\hyperlink{A5l}{A5}), we observe that
    \begin{equation}\label{ken04}
        \frac{1}{2\tau} +g'(\tilde{\eta}) > \frac{1}{4\tau^*} -|g'|_{L^\infty(\R)}  > 0 ~\mbox{ for any $ \tilde{\eta} \in \R $.} 
    \end{equation}
    Because \eqref{ken04} implies the coercivity and strict convexity of the potential $ \Psi_0 $, we conclude the proof of item (I) by applying the general theory of convex analysis (cf. \cite[Proposition 1.2 in Chapter II]{MR1727362}). 

    Additionally, we can verify item (II) immediately as a straightforward consequence of the general theory of convex analysis (cf. \cite[Proposition 1.2 in Chapter II]{MR1727362}) applied to the following proper, l.s.c., coercive, and strictly convex function:
    \begin{equation*}
        \tilde{\theta} \in X \mapsto \frac{1}{2 \tau}|\tilde{\theta} -\bar{\theta}_0|_X^2 +\Phi_{\nu, \varepsilon}(\bar{\eta}_0, \tilde{\theta}) -(\bar{v}, \tilde{\theta})_X \in [0, \infty].
    \end{equation*}

    Finally, we prove item (III). From (\hyperlink{A3l}{A3}),  \eqref{ken03a}, and Remark \ref{uv.P},  we observe that
    \begin{subequations}\label{ken10}
        \begin{align}
            & \frac{1}{\tau} \bigl( (-R_0) -\bar{\eta}_0, {\varphi} \bigr)_X +\kappa^2 \bigl( \nabla (-R_0), \nabla {\varphi} \bigr)_{[X]^N} +\bigl( \alpha'(-R_0) \gamma_\varepsilon(\nabla \bar{\theta}_0), {\varphi} \bigr)_X
            \nonumber
            \\
            & \qquad +\bigl( g(-R_0), \varphi \bigr)_X \leq \bigl( -|\bar{u}|_{L^\infty(\Omega)}, \varphi \bigr)_X,
            \label{ken10a}
            \\[1ex]
            & \frac{1}{\tau} \bigl(R_0 -  \bar{\eta}_0, {\varphi} \bigr)_X +\kappa^2 \bigl( \nabla R_0, \nabla {\varphi} \bigr)_{[X]^N} +\bigl( \alpha'(R_0) \gamma_\varepsilon(\nabla \bar{\theta}_0), {\varphi} \bigr)_X
            \nonumber
            \\
            & \qquad +\bigl( g(R_0), \varphi \bigr)_X \geq \bigl( |\bar{u}|_{L^\infty(\Omega)}, \varphi \bigr)_X 
            \label{ken10b}
            \\[1ex]
            & \mbox{for any $ 0 \leq \varphi \in Y $}
            \nonumber
        \end{align}
    \end{subequations}
    and
    \begin{subequations}\label{ken11}
        \begin{align}
            & \frac{1}{\tau} \bigl( \alpha(\bar{\eta}_0)((-R_0) -\bar{\theta}_0), {\psi} \bigr)_X +\bigl( \alpha(\bar{\eta}_0) \nabla \gamma_\varepsilon (\nabla (-R_0)), \nabla {\psi} \bigr)_{[X]^N} 
            \nonumber
            \\
            & \qquad +M_0 (-R_0) \leq  \bigl( -|\bar{v}|_{L^\infty(\Omega)}, \psi \bigr)_X,
            \label{ken11a}
            \\[1ex]
            & \frac{1}{\tau} \bigl( \alpha(\bar{\eta}_0)(R_0 -\bar{\theta}_0), {\psi} \bigr)_X +\bigl( \alpha(\bar{\eta}_0) \nabla \gamma_\varepsilon (\nabla R_0)), \nabla {\psi} \bigr)_{[X]^N} 
            \nonumber
            \\
            & \qquad +M_0 R_0 \geq  \bigl( |\bar{v}|_{L^\infty(\Omega)}, \psi \bigr)_X
            \label{ken11b}
            \\[1ex]
            & \mbox{for any $ 0 \leq \psi \in Y $.}
            \nonumber
        \end{align}
    \end{subequations}
    We suppose that $ 0 \leq \varphi \in Y $ for the test function in \eqref{ken01} and \eqref{ken10}, and let us take the difference from \eqref{ken10a} to \eqref{ken01}. Then, we observe that
    \begin{align*}
        & \frac{1}{\tau} \bigl( -R_0 -\eta, \varphi \bigr)_X +\kappa^2 \bigl( \nabla (-R_0 -\eta), \nabla \varphi \bigr)_X +\bigl( (\alpha'(-R_0) -\alpha'(\eta)) \gamma_\varepsilon(\nabla \bar{\theta}_0), \varphi \bigr)_X
        \\
        & \qquad +\bigl( g(-R_0) -g(\eta), \varphi \bigr)_X \leq \bigl( -|\bar{u}|_{L^\infty(\Omega)} -u_{i}, \varphi \bigr)_X ~\mbox{ for any $ 0 \leq \varphi \in Y $.}
    \end{align*}
    Subsequently, putting $ \varphi = [-R_0 -\eta]^+ $ in $ X $, we deduce from (\hypertarget{A1l}{A1})--(\hyperlink{A6l}{A6}) and Remark \ref{uv.P} that
    \begin{align}
        \frac{1}{\tau} \bigl| [-R_0 -\eta]^+ \bigr|_X^2 ~& \leq |g'|_{L^\infty(\R)} \bigl| [-R_0 -\eta]^+ \bigr|_X^2 < \frac{1}{4 \tau} \bigl| [-R_0 -\eta]^+ \bigr|_X^2;
        \nonumber
        \\[1ex]
        & \mbox{that is,~ }~ [-R_0 -\eta]^+ = 0 \mbox{ in $ X $.}
        \label{ken12}
    \end{align}
    Simultaneously, taking the difference from \eqref{ken11a} to \eqref{ken02} using a test function $ 0 \leq \psi \in Y $, we compute the following:
    \begin{align*}
        & \frac{1}{\tau} \bigl( \alpha(\bar{\eta}_0) (-R_0 -\theta), \psi \bigr)_X 
        \\
        & +\bigl( \alpha(\bar{\eta}_0) \bigl( \nabla \gamma_\varepsilon(\nabla (-R_0)) -\nabla \gamma_\varepsilon(\nabla \theta) \bigr), \nabla \psi \bigr)_{[X]^N} +\nu^2 \bigl( \nabla \bigl( -R_0 -\theta \bigr), \psi \bigr)_{[X]^N}
        \\
        & +M_0\bigl( -R_0 -\theta, \psi \bigr)_X \leq \bigl( -R_0 -|\bar{v}|_{L^\infty(\Omega)}, \psi \bigr) ~\mbox{ for any $ 0 \leq \psi \in Y $.}
    \end{align*}
    Thus, putting $ \psi = [-R_0 -\theta]^+ $ in $ X $, we infer from (\hypertarget{A1l}{A1})--(\hyperlink{A6l}{A6}) that
    \begin{align}\label{ken13}
        & \delta_*  \bigl| [-R_0 -\theta]^+ \bigr|_X^2 \leq \bigl| \sqrt{\alpha(\bar{\eta}_0)} [-R_0 -\theta]^+ \bigr|_X^2 \leq 0; ~\mbox{ that is, $ [-R_0 -\theta]^+ = 0 $ in $ X $.}
    \end{align}
    \eqref{ken12} and \eqref{ken13} imply that
    \begin{equation}\label{ken14}
        \eta \geq -R_0 ~\mbox{ and }~ \theta \geq -R_0, ~\mbox{ a.e. in $ \Omega $.}
    \end{equation}
    Similarly, we verify the following:
    \begin{align}\label{ken15}
        & [\eta -R_0]^+ = 0  \mbox{ and } [\theta -R_0]^+ = 0  \mbox{ in $ X $;}
        \nonumber
        \\
        & \mbox{that is, }~ \eta \leq R_0 \mbox{ and } \theta \leq R_0 \mbox{, respectively, a.e. in $ \Omega $}
    \end{align}
    by taking the difference from \eqref{ken01} to \eqref{ken10b}, and from \eqref{ken02} to \eqref{ken11b}) under $ \varphi = [\eta -R_0]^+ $ and $ \psi = [\theta -R_0]^+ $, respectively, and using assumptions (\hypertarget{A1l}{A1})--(\hyperlink{A6l}{A6}). 

    \eqref{ken14} and \eqref{ken15} conclude the proof of item (III) because they are equivalent to the estimate \eqref{ken03b}. \qed
}
\begin{lem}\label{Lem3-01}
    We assume (\hyperlink{A1l}{A1})--(\hyperlink{A7l}{A7}). Then, for every $ \nu > 0 $ and $ \varepsilon \in (0, 1) $, the system \hyperlink{(S)$_{\nu, \varepsilon}^{m}$}{(S)$_{\nu, \varepsilon}^{m}$} admits a unique solution $[\{\eta_i \}_{i=1}^m, \{ \theta_i\}_{i=1}^m ] \in [X]^{m \times 2}$, which fulfills the following energy estimate:
    \begin{align}\label{energyEst01}
        \frac{1}{2 \tau} & |\eta_i -\eta_{i -1}|_X^2 +\frac{1}{2 \tau}\bigl| \sqrt{\alpha(\eta_{i -1})} (\theta_i -\theta_{i -1}) \bigr|_{X}^2 +\mathcal{F}_{\nu, \varepsilon}(\eta_i, \theta_i) 
        \nonumber
        \\
        & \leq \mathcal{F}_{\nu, \varepsilon}(\eta_{i-1}, \theta_{i-1}) + \tau R_0^2 \mathcal{L}^N(\Omega) \left( 1 +\frac{1}{2 \delta_*^2} \right), ~\mbox{for i = 1, 2, 3, \dots, m,}
    \end{align}
    where $ R_0 > 0 $ is the constant as in \eqref{ken00} and Remark \ref{uv.P}. Moreover, the system (S)$_{\nu, \varepsilon}^{m} $ admits the following property of $ L^\infty $-bound:
    \begin{itemize}
        \item[$(*)$]if $ [\eta_0, \theta_0] \in [L^\infty(\Omega)]^2 $ and $ |\eta_0|_{L^\infty(\Omega)} \vee |\theta_0|_{L^\infty(\Omega)} \leq R_0 $, then $ \bigl[ \{ \eta_i \}_{i = 1}^m, \{ \theta_i \}_{i = 1}^m \bigr] \in [L^\infty(\Omega)]^{m \times 2} $, and $ |\eta_i|_{L^\infty(\Omega)} \vee |\theta_i|_{L^\infty(\Omega)} \leq R_0 $ for $ i = 1, 2, 3, \dots, m $.
    \end{itemize}
\end{lem}
\vspace{-2ex}
\paragraph{Proof of Lemma \ref{Lem3-01}.}{The existence of solution $ \bigl[ \{\eta_i\}_{i = 1}^m, \{\theta_i\}_{i = 1}^m \bigr] \subset [X]^{m \times 2} $ and the property $(*)$ of $ L^\infty $-bound are immediately verified by applying Lemma \ref{Lem3-02} to the following cases:
\begin{equation*}
    \begin{cases}
        [\bar{\eta}_0, \bar{\theta}_0] = [\eta_{i -1}, \theta_{i -1}] \mbox{ in $ [Y]^2 $,}
        \\
        [\bar{u}, \bar{v}] := [u_i^{m}, v_i^{m}] \mbox{ in $ [L^\infty(\Omega)]^2 $,}
    \end{cases}
\end{equation*}
for every $ i = 1, 2, 3, \dots, m $, inductively. 

Hence, we only have to verify the energy estimate\eqref{energyEst01}.
We fix $i = 1, 2, 3, \ldots, m$ and put $\varphi = \eta_i - \eta_{i-1} $ in \eqref{3tau-1}. Then, we observe that
\begin{align}\label{ken4-20}
    \frac{1}{\tau} |\eta_i - & \eta_{i-1}|_X +\kappa^2 \bigl(\nabla \eta_i, \nabla(\eta_i - \eta_{i-1}) \bigr)_{[X]^N} +\int_\Omega g(\eta_i)(\eta_i - \eta_{i-1}) \, dx \nonumber \\
& +\int_\Omega \alpha'(\eta_i)(\eta_i - \eta_{i-1})\gamma_\varepsilon (\nabla \theta_{i}) \, dx = \bigl(u_i^{m}, \eta_i - \eta_{i-1} \bigr)_X. 
\end{align}
Additionally, by assumption (\hyperlink{A4l}{A4}) and Taylor's theorem, 
\begin{align}\label{ken4-21}
g(\eta_i)(\eta_i - \eta_{i-1}) \geq G(\eta_i)  - G(\eta_{i-1}) - \frac{1}{2}|g'|_{L^\infty(\mathbb{R})}|\eta_i - \eta_{i-1}|^2,  \mbox{ a.e. in } \Omega.
\end{align}
Taking into account \eqref{ken4-21}, the property $(*)$ of $ L^\infty $-bound, assumptions (\hyperlink{A3l}{A3}), (\hyperlink{A4l}{A4}), and (\hyperlink{A5l}{A5}), H\"{o}lder's and Young's inequalities, and Remark \ref{uv.P}, we deduce from \eqref{ken4-20} that
\begin{align}\label{ken4-22}
\displaystyle \frac{1}{2 \tau}|\eta_i - \eta_{i-1}|_X^2 & + \left( \frac{\kappa^2}{2} |\nabla \eta_i|_{[X]^N}^2  - \frac{\kappa^2}{2} |\nabla \eta_{i-1}|_{[X]^N}^2 \right)\nonumber \\
&  +\left( \int_\Omega G(\eta_i)\, dx - \int_\Omega G(\eta_{i-1})\, dx \right) \nonumber \\
& +\left( \int_\Omega \alpha(\eta_i)\gamma_\varepsilon (\nabla \theta_{i})\, dx - \int_\Omega \alpha(\eta_{i-1})\gamma_\varepsilon (\nabla \theta_{i})\, dx  \right) \nonumber \\
& \leq \tau R_0^2 \mathcal{L}^N(\Omega)
\end{align}  
via the following computations:
\begin{align*}
\ds \kappa^2 \bigl( \nabla \eta_i, \nabla(\eta_i -  \eta_{i-1}) \bigr)_{[X]^N} \geq \frac{\kappa^2}{2} |\nabla \eta_i|_{[X]^N}^2  - \frac{\kappa^2}{2} |\nabla \eta_{i-1}|_{[X]^N}^2,
\end{align*}
\begin{align*}
& \int_\Omega g(\eta_i)(\eta_i - \eta_{i-1}) \, dx 
    \\
& \geq \int_\Omega G(\eta_i) \, dx - \int_\Omega G(\eta_{i-1}) \, dx - \frac{1}{2}|g'|_{L^\infty(\mathbb{R})} |\eta_i - \eta_{i-1}|_X^2
\\
& \geq \int_\Omega G(\eta_i) \, dx - \int_\Omega G(\eta_{i-1}) \, dx -\frac{1}{4 \tau} |\eta_i - \eta_{i-1}|_X^2,
\end{align*}
\begin{align*}
& \int_\Omega \alpha'(\eta_i)(\eta_i - \eta_{i-1})\gamma_\varepsilon (\nabla \theta_{i})\, dx 
    \\
& \geq \int_\Omega \alpha(\eta_i)\gamma_\varepsilon (\nabla \theta_{i})\, dx - \int_\Omega \alpha(\eta_{i-1})\gamma_\varepsilon (\nabla \theta_{i})\, dx,
\end{align*}
and
\begin{align*}
    \bigl(u_i^{m}, \eta_i - \eta_{i-1} \bigr)_X & \leq |u_i^{m}|_X |\eta_i - \eta_{i-1}|_X  \leq \tau |u_i^{m}|_X^2 +\frac{1}{4 \tau} |\eta_i -\eta_{i -1}|_X^2
    \\
& \leq \tau R_0^2 \mathcal{L}^N(\Omega) + \frac{1}{4\tau}|\eta_i - \eta_{i-1}|_X^2.
\end{align*}

Next, by putting $\psi = \theta_i - \theta_{i-1} $ in \eqref{3tau-2}, we observe that
\begin{align}\label{ken4-23}
    \ds & \frac{1}{\tau} \bigl| \sqrt{\alpha_0(\eta_{i-1})}(\theta_i - \theta_{i-1}) \bigr|_X^2 
    \nonumber 
    \\
    & \, +\left( \alpha(\eta_{i-1}) \nabla \gamma_\varepsilon(\nabla \theta_i) + \nu^2 \nabla \theta_i, \nabla (\theta_i - \theta_{i-1}) \right)_{[X]^N} 
    \nonumber 
    \\
    & \, + M_0 \bigl( \theta_i, \theta_i - \theta_{i-1} \bigr)_X = \bigl(v_i^{m}, \theta_i - \theta_{i-1} \bigr)_X.
\end{align}
Thus, bearing in mind the property $(*)$ of $L^\infty$-bound, assumption (\hyperlink{A3l}{A3}), and Remark \ref{uv.P}, and using Young's inequality, we infer from \eqref{ken4-23} that
 \begin{align}\label{ken4-24}
\ds & \frac{1}{2\tau}|\sqrt{\alpha_0(\eta_{i-1})}(\theta_i - \theta_{i-1})|_X^2 \nonumber \\
&\, +\left(\int_\Omega \alpha(\eta_{i-1}) \gamma_\varepsilon(\nabla \theta_i) \, dx - \int_\Omega \alpha(\eta_{i-1}) \gamma_\varepsilon(\nabla \theta_{i-1})\, dx \right)
     \nonumber 
     \\
     &\, +\frac{\nu^2}{2} \bigl( |\nabla \theta_i| _{[X]^N}^2 - |\nabla \theta_{i-1}|_{[X]^N}^2 \bigr) + \frac{M_0}{2} \bigl( |\theta_i|_X^2 - |\theta_{i-1}|_X^2 \bigr) 
     \leq \frac{\tau R_0^2 \mathcal{L}^N(\Omega)}{2 \delta_*^2}
\end{align}
via the following computations:
\begin{align*}
\ds & \left( \alpha(\eta_{i-1}) \nabla \gamma_\varepsilon(\nabla \theta_i) + \nu^2 \nabla \theta_i, \nabla (\theta_i - \theta_{i-1}) \right)_{[X]^N} 
    \\
    & = \int_\Omega \alpha(\eta_{i-1}) \nabla \gamma_\varepsilon(\nabla \theta_i) \cdot  \nabla (\theta_i - \theta_{i-1})\, dx +\nu^2 \bigl(\nabla \theta_i,  \nabla (\theta_i - \theta_{i-1})\bigr)_{[X]^N} 
    \\
    & \geq \int_\Omega \alpha(\eta_{i-1})\gamma_\varepsilon(\nabla \theta_i) \, dx - \int_\Omega \alpha(\eta_{i-1}) \gamma_\varepsilon(\nabla \theta_{i-1})\, dx 
    \\
   & \quad + \frac{\nu^2}{2} \bigl(|\nabla \theta_i| _{[X]^N}^2 - |\nabla \theta_{i-1}|_{[X]^N}^2 \bigr),
\end{align*}
\begin{align*}
\ds M_0 \bigl(\theta_i, \theta_i - \theta_{i-1} \bigr)_X \geq \frac{M_0}{2} \bigl( |\theta_i|_X^2 - |\theta_{i-1}|_X^2 \bigr),
\end{align*}
and
\begin{align*}
    \ds \bigl( v_i^{m}, \theta_i - \theta_{i-1} \bigr)_X & \leq |v_i^{m}|_X \cdot \frac{1}{\delta_*}\bigl| \sqrt{\alpha_0(\eta_{i -1})} (\theta_i -\theta_{i -1}) \bigr|_X 
    \\
    & \leq \frac{\tau}{2\delta_*^2} |v_i^{m}|_X^2 + \frac{1}{2 \tau} \bigl| \sqrt{\alpha_0}(\eta_{i -1}) (\theta_i -\theta_{i -1}) \bigr|_X^2 
    \\
    & \leq \frac{\tau R_0^2 \mathcal{L}^N(\Omega) }{2\delta_*^2} +\frac{1}{2 \tau} \bigl| \sqrt{\alpha_0}(\eta_{i -1}) (\theta_i -\theta_{i -1}) \bigr|_X^2.
\end{align*}
We verify the energy estimate\eqref{energyEst01} by summing \eqref{ken4-22} and \eqref{ken4-24}.

Thus, we complete the proof of the lemma. 
\qed
}

\begin{lem}\label{Lem3-03KS}
We assume  (\hyperlink{A1l}{A1})--(\hyperlink{A7l}{A7}), and $\nu \in (0, \nu_0 + 1]$. We consider
    the solution $[\{\eta_i \}_{i=1}^m, \{\theta_i \}_{i=1}^m] \in  [X]^{m \times 2} $ to the approximating system (S)$_{\nu, \varepsilon}^{m}$. 
    Then, there exist constants $ R_1, R_2, R_* > 0$, independent of $\nu$, $ \varepsilon $, $ m $, and the initial pair $ [\eta_0, \theta_0] $, such that a sequence $ \{ \mathcal{X}_i \}_{i = 1}^m \subset [0, \infty) $, defined as
\begin{align}\label{risan.gron}
    \ds \mathcal{X}_i := |\eta_i|_X^2 + R_* |\theta_i|_X^2 +  \kappa^2\tau|\nabla \eta_{i}|_{[X]^N}^2, ~ i = 0, 1, 2, 3, \dots, m,
\end{align}

    fulfills the following properties: 
 \begin{description}
     \item[\textmd{(\hypertarget{a}{a})}]$\ds \mathcal{X}_0 \leq R_1 $ implies $ \mathcal{X}_i \leq R_1 $ for all $ i = 1, 2, 3, \dots, m $, and in particular,
         \begin{align*}
             & |\eta_m|_X^2 + R_* |\theta_m|_X^2 \leq \mathcal{X}_m \leq R_1. 
         \end{align*} 
     \item[\textmd{(\hypertarget{b}{b})}]$\ds \mathcal{X}_0  \leq R_1$ implies $ i \tau \mathcal{F}_{\nu, \varepsilon}(\eta_i, \theta_i) \leq R_2$ for all $ i = 1, 2, 3, \dots, m $, and in particular, 
         \begin{align*}
             &
             \kappa^2 |\nabla \eta_m|_{[X]^N}^2 +\delta_* |D \theta_m|(\Omega) +\nu^2 |\nabla \theta_m|_{[X]^N}^2 \leq \frac{2 R_2}{T} := R_3.
         \end{align*}
\end{description}
\end{lem}

\paragraph{Proof of Lemma \ref{Lem3-03KS}.}{
        We fix $i \in \{ 1, 2, 3, \dots, m \} $ and fix the positive constant $L > 0$. Additionally, we put $\varphi = \eta_i$ in \eqref{3tau-1} and add $(L \eta_i, \eta_i)_X$ to both sides of the result. Then, we obtain 
\begin{align}\label{4-03}
    \displaystyle \frac{1}{\tau} & (\eta_i - \eta_{i-1}, \eta_i)_X + \kappa^2|\nabla \eta_i|_{[X]^N}^2 + \int_\Omega g(\eta_i)\eta_i \, dx \nonumber \\
& + \int_\Omega \alpha'(\eta_i)\gamma_\varepsilon (\nabla \theta_{i})\eta_i\, dx + L|\eta_i|_X^2 = (u_i^{m}, \eta_i)_X + L(\eta_i, \eta_i)_X. 
\end{align}
Additionally, given (\hyperlink{A4l}{A4}) and Taylor's theorem, 
\begin{align}\label{4-04}\noeqref{4-04}
g(\eta_i)\eta_i \geq G(\eta_i)  - G(0) - \frac{1}{2}|g'|_{L^\infty(\mathbb{R})}|\eta_i|^2,  \mbox{ a.e. in } \Omega.
\end{align}

Note that we set the positive constant $C_1$, independent of $\nu, \varepsilon, \mbox{ and } m$, as follows: 
\begin{align}\label{4-04-1}
 C_1 :=  |G(0)|_{L^1(\Omega)} + R_0^2 \mathcal{L}^N(\Omega) + \frac{R_0^2 \mathcal{L}^N(\Omega)}{2}|g'|_{L^\infty(\mathbb{R})} + L R_0^2 \mathcal{L}^N(\Omega).
\end{align} 
Using assumptions (\hyperlink{A3l}{A3})--(\hyperlink{A5l}{A5}), Remark \ref{uv.P}, Lemma \ref{Lem3-01}, and Young's inequality, we deduce from \eqref{4-03}--\eqref{4-04-1} that
\begin{align}\label{4-05}
\displaystyle & \frac{1}{2\tau}(|\eta_i|_X^2 - |\eta_{i-1}|_X^2) + \kappa^2|\nabla \eta_i|_{[X]^N}^2 + \int_\Omega G(\eta_i)\, dx + L|\eta_i|_X^2 \leq C_1 
\end{align} 
via
\begin{align*}
\ds \frac{1}{\tau} (\eta_i - \eta_{i-1}, \eta_i)_X \geq \frac{1}{2\tau}(|\eta_i|_X^2 - |\eta_{i-1}|_X^2), 
\end{align*}
\begin{align*}
\int_\Omega g(\eta_i)\eta_i \, dx & \geq \int_\Omega G(\eta_i) \, dx - \int_\Omega G(0) \, dx - \frac{1}{2}|g'|_{L^\infty(\mathbb{R})}\int_\Omega |\eta_i|^2 \, dx \\
& \geq \int_\Omega G(\eta_i) \, dx - \int_\Omega G(0) \, dx - \frac{R_0^2 \mathcal{L}^N (\Omega) }{2}|g'|_{L^\infty(\mathbb{R})},
\end{align*}
\begin{align*}
\int_\Omega \alpha'(\eta_i)\gamma_\varepsilon (\nabla \theta_{i})\eta_i\, dx = \int_\Omega (\alpha'(\eta_i) - \alpha'(0))(\eta_i - 0)\gamma_\varepsilon (\nabla \theta_{i})\, dx \geq 0,
\end{align*}
\begin{align*}
\ds (u_i^{m}, \eta_i)_X = \int_\Omega u_i^{m} \eta_i\, dx \leq R_0^2 \int_\Omega \, dx = R_0^2 \mathcal{L}^N(\Omega),
\end{align*}
and
\begin{align*}
\ds L(\eta_i, \eta_i)_X = L \int_\Omega|\eta_i|^2\, dx \leq L R_0^2 \mathcal{L}^N(\Omega).
\end{align*}

Next, we consider $\psi = \theta_i/\alpha_0(\eta_{i-1})$ in \eqref{3tau-2} and compute the following:
\begin{align}\label{4-06}
\ds \frac{1}{\tau}(\theta_i - \theta_{i-1}, \theta_i)_X & + \left( \alpha(\eta_{i-1}) \nabla \gamma_\varepsilon(\nabla \theta_i) + \nu^2 \nabla \theta_i, \nabla \left( \frac{\theta_i}{\alpha_0(\eta_{i-1})}\right) \right)_{[X]^N} \nonumber \\
& + \left(M_0\theta_i, \frac{\theta_i}{\alpha_0(\eta_{i-1})} \right)_X = \left(v_i^{m}, \frac{\theta_i}{\alpha_0(\eta_{i-1})} \right)_X.
\end{align}

Taking into account assumption (\hyperlink{A3l}{A3}), Example \ref{Rem.gamma_eps}, Remark \ref{uv.P}, Lemma \ref{Lem3-01}, and Young's inequality, we estimate from \eqref{4-06} that
 \begin{align}\label{4-07}
     & \ds \frac{1}{2\tau}(|\theta_i|_X^2 - |\theta_{i-1}|_X^2) + \frac{\delta_*}{|\alpha_0|_{L^\infty(-R_0, R_0)} |\alpha|_{L^\infty(-R_0, R_0)}  } \int_\Omega \alpha(\eta_{i}) |\nabla \theta_i|\, dx  \nonumber \\
& \quad + \frac{\nu^2}{2|\alpha_0|_{L^\infty(-R_0, R_0)}}|\nabla \theta_i|_{[X]^N}^2 + \frac{M_0}{|\alpha_0|_{L^\infty(-R_0, R_0)}}|\theta_i|_X^2 \nonumber \\
 \leq~ & \frac{R_0^2}{\delta_*}\mathcal{L}^N(\Omega) + \frac{\kappa^2}{2}|\nabla \eta_{i-1}|_{[X]^N}^2  \nonumber \\
& \quad + \frac{1}{2\kappa^2}\left( \frac{ |\alpha|_{L^\infty(-R_0, R_0)} |\alpha_0'|_{L^\infty(-R_0, R_0)}}{\delta_{*}^{2}} R_0 \right)^2 \mathcal{L}^N(\Omega) \nonumber \\
& \quad + \frac{\nu^2 |\alpha_0|_{L^\infty(-R_0, R_0)}}{2} \left( \frac{|\alpha_0'|_{L^\infty(-R_0, R_0)}^2 R_0}{\delta_*^2 }\right)^2 |\nabla \eta_{i-1}|_{[X]^N}^2
\end{align}
via
\begin{align*}
\ds \frac{1}{\tau}(\theta_i - \theta_{i-1}, \theta_i)_X \geq \frac{1}{2\tau}(|\theta_i|_X^2 - |\theta_{i-1}|_X^2)
\end{align*}
\begin{align*}
& \ds \left( \alpha(\eta_{i-1}) \nabla \gamma_\varepsilon(\nabla \theta_i) + \nu^2 \nabla \theta_i, \nabla \left( \frac{\theta_i}{\alpha_0(\eta_{i-1})}\right) \right)_{[X]^N} \\
& \qquad = \left( \frac{\alpha(\eta_{i-1})}{\alpha_0(\eta_{i-1})}  \nabla \gamma_\varepsilon(\nabla \theta_i), \nabla \theta_i\right)_{[X]^N} + \left( \frac{\nu^2 \nabla \theta_i}{\alpha_0(\eta_{i-1})}, \nabla \theta_i\right)_{[X]^N} \\
& \qquad \qquad  - \left( \frac{\alpha(\eta_{i-1}) \alpha_0'(\eta_{i-1}) \theta_i}{\alpha_0(\eta_{i-1})^2}  \nabla \gamma_\varepsilon(\nabla \theta_i), \nabla \eta_{i-1}\right)_{[X]^N} 
    \\
    & \qquad \qquad 
    - \left( \frac{\nu^2 \alpha_0'(\eta_{i-1}) \theta_i}{\alpha_0(\eta_{i-1})^2} \nabla \eta_{i-1},  \nabla \theta_{i} \right)_{[X]^N},
\end{align*}
\begin{align*}
    \ds \left( \frac{\alpha(\eta_{i-1})}{\alpha_0(\eta_{i-1})} \right. & \left. \rule{0pt}{16pt} \nabla \gamma_\varepsilon(\nabla \theta_i), \nabla \theta_i \right)_{[X]^N} = \int_\Omega \frac{\alpha(\eta_{i-1})}{\alpha_0(\eta_{i-1})}  \nabla \gamma_\varepsilon(\nabla \theta_i) \cdot \nabla \theta_i \, dx \\
& \geq \frac{1}{|\alpha_0|_{L^\infty(-R_0, R_0)} }\int_\Omega \alpha(\eta_{i-1})|\nabla \theta_i|\, dx \\
& \geq  \frac{\delta_*}{|\alpha_0|_{L^\infty(-R_0, R_0)} |\alpha|_{L^\infty(-R_0, R_0)}  } \int_\Omega \alpha(\eta_{i}) |\nabla \theta_i|\, dx,
\end{align*}
\begin{align*}
\ds \left( \frac{\nu^2 \nabla \theta_i}{\alpha_0(\eta_{i-1})}, \nabla \theta_i\right)_{[X]^N} \geq \frac{\nu^2}{|\alpha_0|_{L^\infty(-R_0, R_0)}} |\nabla \theta_i|_{[X]^N}^2, 
\end{align*}
\begin{align*}
    \ds - & \left( \frac{\alpha(\eta_{i-1}) \alpha_0'(\eta_{i-1}) \theta_i}{\alpha_0(\eta_{i-1})^2} \nabla \gamma_\varepsilon(\nabla \theta_i), \nabla \eta_{i-1} \right)_{[X]^N} 
    \\
    & \geq - \int_\Omega \frac{ |\alpha|_{L^\infty(-R_0, R_0)} |\alpha_0'|_{L^\infty(-R_0, R_0)}}{\delta_{*}^{2}} R_0 |\nabla \eta_{i-1}|\, dx \\
& \geq - \frac{\kappa^2}{2}|\nabla \eta_{i-1}|_{[X]^N}^2 - \frac{1}{2\kappa^2}\left( \frac{ |\alpha|_{L^\infty(-R_0, R_0)} |\alpha_0'|_{L^\infty(-R_0, R_0)}}{\delta_{*}^{2}} R_0 \right)^2 \mathcal{L}^N(\Omega),
\end{align*}
\begin{align*}
    \ds - & \left( \frac{\nu^2 \alpha_0'(\eta_{i-1}) \theta_i}{\alpha_0(\eta_{i-1})^2} \nabla \eta_{i-1},  \nabla \theta_i\right)_{[X]^N} \geq - \int_\Omega \frac{ \nu^2 |\alpha_0'|_{L^\infty(-R_0, R_0)}}{\delta_{*}^{2}} R_0 |\nabla \theta_i| |\nabla \eta_{i-1}|\, dx \\ 
& \geq - \frac{\nu^2}{2|\alpha_0|_{L^\infty(-R_0, R_0)}}|\nabla \theta_i|_{[X]^N}^2 
    - \frac{\nu^2 |\alpha_0|_{L^\infty(-R_0, R_0)}}{2} \left( \frac{|\alpha_0'|_{L^\infty(-R_0, R_0)} R_0}{\delta_*^2 }\right)^2 |\nabla \eta_{i-1}|_{[X]^N}^2,
\end{align*}
\begin{align*}
\ds \left(M_0\theta_i, \frac{\theta_i}{\alpha_0(\eta_{i-1})} \right)_X \geq \frac{M_0}{|\alpha_0|_{L^\infty(-R_0, R_0)}}|\theta_i|_X^2,
\end{align*}
and
\begin{align*}
\ds \left(v_i^{m}, \frac{\theta_i}{\alpha_0(\eta_{i-1})} \right)_X & \leq \int_\Omega \frac{1}{\delta_*}|v_i^{m}||\theta_i|\, dx \leq \frac{R_0^2}{\delta_*}\mathcal{L}^N(\Omega).
\end{align*}
Based on this, we define the positive constant $R_* \in (0, 1)$ by putting
\begin{align}\label{4-08}
\ds R_* := & \frac{\kappa^2}{2(1 + \kappa^2)} \times \frac{1}{1 + R_0^2} \times \frac{\delta_*^4}{1 + \delta_*^4} \nonumber \\
 & \times \frac{1}{2(1 +\nu_0)^2  }\times \frac{1}{1 + |\alpha_0|_{L^\infty(-R_0, R_0)}} \times \frac{1}{1 + |\alpha_0'|_{L^\infty(-R_0, R_0)}^{2}}.
\end{align}
Additionally, we multiply both sides of \eqref{4-07} by $R_*$. Then, we observe that
\begin{align}\label{4-09}
\ds & \frac{R_*}{2\tau}(|\theta_i|_X^2 - |\theta_{i-1}|_X^2) +  \frac{\delta_* R_*}{|\alpha_0|_{L^\infty(-R_0, R_0)} |\alpha|_{L^\infty(-R_0, R_0)}  } \int_\Omega \alpha(\eta_{i}) |\nabla \theta_i|\, dx  \nonumber \\
& \quad + \frac{\nu^2 R_*}{2|\alpha_0|_{L^\infty(-R_0, R_0)}}|\nabla \theta_i|_{[X]^N}^2 + \frac{M_0 R_*}{|\alpha_0|_{L^\infty(-R_0, R_0)}}|\theta_i|_X^2 \nonumber \\
 & \leq R_* \mathcal{L}^N(\Omega) \left( \frac{R_0^2}{\delta_*} + \frac{1}{2\kappa^2}\left( \frac{ |\alpha|_{L^\infty(-R_0, R_0)} |\alpha_0'|_{L^\infty(-R_0, R_0)}}{\delta_{*}^{2}} R_0 \right)^2  \right) + \frac{\kappa^2}{2}|\nabla \eta_{i-1}|_{[X]^N}^2 \nonumber \\
& = C_2 + \frac{\kappa^2}{2}|\nabla \eta_{i-1}|_{[X]^N}^2,
\end{align}
where
\begin{align*}
 C_2 :=  R_* \mathcal{L}^N(\Omega) \left( \frac{R_0^2}{\delta_*} + \frac{1}{2\kappa^2}\left( \frac{ |\alpha|_{L^\infty(-R_0, R_0)} |\alpha_0'|_{L^\infty(-R_0, R_0)}}{\delta_{*}^{2}} R_0 \right)^2  \right). 
\end{align*}
Note that $C_2$ is independent of $\nu $, $ \varepsilon $, and $m$. 

By summing the inequalities \eqref{4-05} and \eqref{4-09}, and invoking \eqref{risan.gron}, we obtain \begin{align}\label{4-10}
    \ds & \frac{1}{2\tau} \bigl( \mathcal{X}_i -\mathcal{X}_{i -1} \bigr) + L|\eta_i|_X^2 + \frac{\kappa^2}{2}|\nabla \eta_i|_{[X]^N}^2 + \int_\Omega G(\eta_i)\, dx \nonumber 
    \\[1ex]
& \quad + \frac{\delta_* R_*}{|\alpha_0|_{L^\infty(-R_0, R_0)} |\alpha|_{L^\infty(-R_0, R_0)} } \int_\Omega \alpha(\eta_i) |\nabla \theta_i|\, dx + \frac{\nu^2 R_*}{2|\alpha_0|_{L^\infty(-R_0, R_0)}}|\nabla \theta_i|_{[X]^N}^2  \nonumber 
    \\[1ex]
& \quad + \frac{M_0 R_*}{|\alpha_0|_{L^\infty(-R_0, R_0)}}|\theta_i|_X^2 + \int_\Omega G(\eta_i) \, dx  \leq C_1 + C_2  =: C_3.
\end{align}
Additionally, we define the positive constant $C_4 \in (0, 1]$ as follows: 
\begin{align}\label{4-11}
\ds C_4 :=  \frac{\delta_* }{|\alpha_0|_{L^\infty(-R_0, R_0)} |\alpha|_{L^\infty(-R_0, R_0)} }  \wedge \frac{1}{ |\alpha_0 |_{L^\infty(-R_0, R_0)}}.
\end{align}
Note that $C_4$ is independent of $\nu $, $ \varepsilon $, and $ m$. 
From (\hyperlink{A6l}{A6}), \eqref{4-10}, and \eqref{4-11}, it follows that
\begin{align}\label{4-11-1}
 \ds & \frac{1}{2\tau} \bigl( \mathcal{X}_i -\mathcal{X}_{i -1} \bigr) + L|\eta_i|_X^2 + \frac{\kappa^2}{2} |\nabla \eta_i|_{[X]^N}^2 + \int_\Omega G(\eta_i)\, dx + C_4 R_* \Phi_{\nu, \varepsilon}(\eta_i, \theta_i)  \leq C_3. 
\end{align}

Now, we verify item  (\hyperlink{a}{a}). 
From assumptions (\hyperlink{A4l}{A4})--(\hyperlink{A6l}{A6}), and \eqref{4-11-1}, we observe that 
\begin{align}\label{4-12}
    \ds & \frac{1}{2\tau} \bigl( \mathcal{X}_i -\mathcal{X}_{i -1} \bigr) + L|\eta_i|_X^2 + \frac{\kappa^2 \tau}{2} |\nabla \eta_i|_{[X]^N}^2 + \frac{C_4 R_* M_0}{2}|\theta_i|_X^2 \leq C_3. 
\end{align}
We set the positive constant $C_5 \in (0, 1]$ as follows: 
\begin{align}\label{4-12-1}
\ds C_5 := \frac{1}{2} \wedge L \wedge \frac{C_4 M_0}{2}.
\end{align}
Note that $C_5$ is independent of $\nu$, $\varepsilon$, and $ m$. 
Taking into account \eqref{risan.gron}, \eqref{4-12}, and \eqref{4-12-1}, we estimate the following: 
\begin{align}\label{4-12-2}
    \ds & \frac{1}{2\tau} \bigl( \mathcal{X}_i -\mathcal{X}_{i -1} \bigr) + C_5 \mathcal{X}_i \leq C_3. 
\end{align}

Now, applying Lemma \ref{Ap.Lem1} to the case when 
\begin{align*}
    \ds \{ A_i \} = \{ \mathcal{X}_i \}, ~ \Lambda_* = C_5, ~\mbox{and}~ K_* = C_3,
\end{align*}
we deduce from \eqref{4-12-2} that
        \begin{align}\label{4-13}
            \mathcal{X}_i \leq \left( \frac{1}{1 + 2\tau C_5} \right)^i \mathcal{X}_0 + \frac{C_3}{C_5}\left[ 1 - \left( \frac{1}{1 + 2\tau C_5} \right)^i \right], \mbox{ for  } i = 1, 2, 3, \ldots, m.  
        \end{align}
Thus, by setting
\begin{align}\label{4-14}
 \ds R_1 := \frac{C_3}{C_5}, \mbox{ independent of } \nu, \varepsilon \mbox{ and } m,
\end{align}
we verify item (\hyperlink{a}{a}) as a consequence of \eqref{4-13}.
\medskip

Next, to verify (\hyperlink{b}{b}), we define the positive constant $C_6 \in (0, 1]$ as follows: 
\begin{align}\label{4-14-1} 
 \ds C_6 := 1 \wedge C_4 R_*. 
\end{align}
Note that $C_6$ is independent of $\nu$, $\varepsilon$, and  $ m$. 
Using \eqref{free.en}, \eqref{4-11-1}, and \eqref{4-14-1}, we deduce the following: 
\begin{align}\label{4-14-2}
 \ds \frac{1}{2\tau} \bigl( \mathcal{X}_i -\mathcal{X}_{i -1} \bigr) + C_6 \mathcal{F}_{\nu, \varepsilon}(\eta_i, \theta_i) \leq C_3.  
\end{align}
Now, we multiply both sides of \eqref{4-14-2} by $ \tau/C_6 $ and sum the inequalities \eqref{4-14-2} for $i = 1, 2, 3, \ldots, m$. 
Then, bearing in mind \eqref{risan.gron} and (\hyperlink{A5l}{A5}), we compute the following: 
\begin{align*}
    \ds & \frac{1}{2C_6} \bigl( \mathcal{X}_i -\mathcal{X}_{0} \bigr) + \tau \sum_{i=1}^m \mathcal{F}_{\nu, \varepsilon}(\eta_i, \theta_i)  \leq m \cdot \frac{\tau C_3}{C_6} = \frac{T C_3}{C_6}; 
\end{align*}
therefore, 
\begin{align}\label{4-16}
    \tau \sum_{i = 1}^m & \mathcal{F}_{\nu, \varepsilon}(\eta_i, \theta_i)  \leq \frac{1}{2C_6} \mathcal{X}_0 + \frac{T C_3}{C_6}.
\end{align}
Taking into account (\hyperlink{a}{a}), \eqref{4-14}, and \eqref{4-16}, we estimate the following: 
\begin{align}\label{4-17}
    \tau \sum_{i = 1}^m \mathcal{F}_{\nu, \varepsilon}(\eta_i, \theta_i) \leq ~& \frac{R_1}{2C_6} + \frac{T C_3}{C_6} =: C_7
    \\
    & \mbox{whenever $ \mathcal{X}_0 \leq R_1 $.}
    \nonumber
\end{align}
Furthermore, by multiplying both sides of \eqref{energyEst01} by $ (i -1) \tau $, we deduce the following:
\begin{align*}
 \ds & \frac{i-1}{2} \left(|\eta_i - \eta_{i-1}|_X^2 +  |\sqrt{\alpha_0(\eta_{i-1})}(\theta_i - \theta_{i-1})|_X^2 \right) + i\tau \mathcal{F}_{\nu, \varepsilon}  (\eta_i, \theta_i) \nonumber \\
      & \leq (i-1)\tau\mathcal{F}_{\nu, \varepsilon}  (\eta_{i-1}, \theta_{i-1}) + (i-1)\tau^2  R_0^2 \mathcal{L}^N(\Omega) \left( 1 +\frac{1}{2 \delta_*^2} \right)  \nonumber \\ 
 & \quad + \tau \mathcal{F}_{\nu, \varepsilon}  (\eta_i, \theta_i), \mbox{ for } i = 1, 2, 3, \ldots, m,  
\end{align*}
which implies that
\begin{align}\label{4-26}
    \ds i\tau \mathcal{F}_{\nu, \varepsilon}  (\eta_i, \theta_i) & \leq \sum_{j = 1}^i (j -1) \tau^2  R_0^2 \mathcal{L}^N(\Omega) \left( 1 +\frac{1}{2 \delta_*^2} \right)  + \tau \sum_{j=1}^i  \mathcal{F}_{\nu, \varepsilon}  (\eta_j, \theta_j) \nonumber \\
& \leq \sum_{i = 1}^m (i -1) \tau^2  R_0^2 \mathcal{L}^N(\Omega) \left( 1 +\frac{1}{2 \delta_*^2} \right)  + \tau \sum_{i=1}^m  \mathcal{F}_{\nu, \varepsilon}  (\eta_i, \theta_i), \\
 & \qquad \mbox{ for } i = 1, 2, 3, \ldots, m. \nonumber
\end{align} 

Now, putting
\begin{align}\label{4-28}
 R_2 := & T^2 R_0^2 \mathcal{L}^N(\Omega) \left( 1 +\frac{1}{2 \delta_*^2} \right)  +C_7,
\end{align}
we verify item (b) as a consequence of \eqref{4-17} and \eqref{4-26}. 
 
Thus, we conclude the proof of this lemma. 
\qed

}

\begin{lem}\label{Lem-Gamma01}
    Under assumptions (\hyperlink{A1l}{A1})--(\hyperlink{A7l}{A7}), let $ \bar{\nu}  \in (0, \nu_0 + 1]$ and $ \bar{\varepsilon} \in (0, 1) $ be fixed constants. Let $ \{ \bar{\nu}_n \}_{n = 1}^\infty \subset (0, \nu_0 + 1] $ and $ \{ \bar{\varepsilon}_n \}_{n = 1}^\infty \subset (0, 1) $ be sequences such that
    \begin{align}\label{bar-nu-eps}
        & \bar{\nu}_n \to \bar{\nu} \mbox{ and } \bar{\varepsilon}_n \to \bar{\varepsilon} \mbox{ as $ n \to \infty $.}
    \end{align}
    Additionally, we assume  
    $ \bar{\eta} \in L^\infty(\Omega) $, $ \{ \bar{\eta}_n \}_{n = 1}^\infty \subset L^\infty(\Omega) $, and
    \begin{equation}\label{Gconv01}
    \left\{ \hspace{-3ex} \parbox{8.5cm}{
        \vspace{-2ex}
        \begin{itemize}
        \item$ \{ \bar{\eta}_{n}\}_{n = 1}^\infty  $ is bounded in $ L^\infty(\Omega) $, 
            \vspace{-1ex}
            \item$ \bar{\eta}_{n} \to \bar{\eta} $ in the pointwise sense, a.e. in $ \Omega $, as $ n  \to \infty $. 
        \end{itemize}
        \vspace{-2ex}
    } \right. 
\end{equation}
    Then, for the sequence of convex functions $ \{ \Phi_{\bar{\nu}_n, \bar{\varepsilon}_n}(\bar{\eta}_n, \cdot\,) \}_{n = 1}^\infty $, it holds that $ \Phi_{\bar{\nu}_n, \bar{\varepsilon}_n}(\bar{\eta}_n, \cdot\,) \to  \Phi_{\nu, \varepsilon}(\bar{\eta}, \cdot\,) $ on $ X $, in the sense of Mosco, as $ n \to \infty $. 
            Moreover, for any $ \bar{\theta} \in D(\Phi_{\bar{\nu}, \bar{\varepsilon}}(\bar{\eta}, \cdot)) $ $ (= L^2(I; Y)) $, and a sequence $ \{ \bar{\theta}_n \}_{n = 1}^\infty \subset X $, the convergence of optimality
            \begin{align}\label{Gconv02}
                & \Phi_{\bar{\nu}_n, \bar{\varepsilon}_n}(\bar{\eta}_n, \bar{\theta}_n) \to \Phi_{\bar{\nu}, \bar{\varepsilon}}(\bar{\eta}, \bar{\theta}) \mbox{ as $ n \to \infty $} \end{align}
        implies that
            \begin{align*}
                & \bar{\theta}_n \to \bar{\theta} \mbox{ in $ Y $ as $ n \to \infty $.}
        \end{align*}
\end{lem}

    \paragraph{Proof of Lemma \ref{Lem-Gamma01}.}{First, we prove the Mosco convergence of the sequence of convex functions $ \{ \Phi_{\bar{\nu}_n, \bar{\varepsilon}_n}(\bar{\eta}_n,\cdot{}) \}_{n = 1}^\infty $.

    To prove the condition of the lower bound, we consider a sequence of functions $ \{ \check{\theta}_n \}_{n = 1}^\infty \subset X $ with a weak limit  $ \check{\theta} $ in $ X $; that is,
    \begin{align*}
        & \check{\theta}_n \to \check{\theta} \mbox{ weakly in $ X $ as $ n \to \infty $.}
    \end{align*}
    Additionally, we consider a subsequence $ \{n_k\}_{k = 1}^\infty \subset \{n\} $ such that
    \begin{align*}
        & \Lambda_0 := \liminf_{n \to \infty} \Phi_{\bar{\nu}_n, \bar{\varepsilon}_n}(\bar{\eta}_n, \check{\theta}_n) = \lim_{k \to \infty} \Phi_{\bar{\nu}_{n_k}, \bar{\varepsilon}_{n_k}}(\bar{\eta}_{n_k}, \check{\theta}_{n_k}) \in [0, \infty].
    \end{align*}
    Then, it is sufficient to consider only the case when $ \Lambda_0 < \infty $ because any other case is obvious. In this case, assumptions \eqref{bar-nu-eps}, (\hyperlink{A3l}{A3}), and (\hyperlink{A6l}{A6}) guarantee the boundedness of $ \{ \check{\theta} \}_{n = 1}^\infty  $ in $ Y $. Hence, bearing \eqref{bar-nu-eps} and \eqref{Gconv01} in mind, we infer that
    \begin{align}\label{KenM01}
            \check{\theta}_{n_k} \to \check{\theta} \mbox{ in $ X $, weakly in $ Y $, and in the pointwise sense, a.e. in $ \Omega $,} \mbox{ as $ k \to \infty $}
    \end{align}
    by considering subsequence(s) again if necessary. The convergence \eqref{KenM01} and the (weakly) lower semi-continuity of the $L^1$-norm and $ L^2 $-based norm on $ X $ lead to
\begin{align}\label{*ken}
\ds \liminf_{k \to \infty} &\,  \bigl| {\ts \alpha(\bar{\eta}_{n_k})} \gamma_{\bar{\varepsilon}_{n_k}}(\nabla \check{\theta}_{n_k}) \bigr|_{L^1(\Omega)} \nonumber \\ 
& \geq \liminf_{k \to \infty} \int_\Omega \alpha (\bar{\eta}) \gamma_{\bar{\varepsilon}}(\nabla \check{\theta}_{n_k})\, dx  + \liminf_{k \to \infty} \Bigl(- \int_\Omega \alpha(\bar{\eta}_{n_k})\bigl|\gamma_{\bar{\varepsilon_{n_k}}}(\nabla \check{\theta}_{n_k} )- \gamma_{\bar{\varepsilon}}(\nabla \check{\theta}_{n_k}) \bigr| \, dx \nonumber \\
& \quad \quad \quad - \int_\Omega \bigl|\alpha(\bar{\eta}_{n_k}) -  \alpha(\bar{\eta}) \bigr| \gamma_{\bar{\varepsilon}}(\nabla \check{\theta}_{n_k})\, dx  \Bigr) \nonumber \\
& \geq \int_\Omega \alpha (\bar{\eta}) \gamma_{\bar{\varepsilon}}(\nabla \check{\theta})\, dx  - \left( \lim_{k \to \infty} |\bar{\varepsilon}_{n_k} - \bar{\varepsilon}|_{L^1(\Omega)} \right) \sup_{n \in \mathbb{N}}|\alpha (\bar{\eta}_n)|_{L^\infty(\Omega)} \nonumber \\
& \, \quad - |\alpha'|_{L^\infty(-R_0, R_0)} \left(\lim_{k \to \infty} |\bar{\eta}_{n_k} - \bar{\eta}|_X \right) \sup_{n \in \mathbb{N}} \bigl|\bar{\varepsilon} + |\nabla \check{\theta}_n|\bigr|_X \nonumber \\
& = \bigl| \alpha(\bar{\eta}) \gamma_{\bar{\varepsilon}}(\nabla \check{\theta})\bigr|_{L^1(\Omega)}.
\end{align}
Hence,
    \begin{align*}
        & \liminf_{n \to \infty} \Phi_{\bar{\nu}_n, \bar{\varepsilon}_n}(\bar{\eta}_n, \check{\theta}_n) ~ (= \Lambda_0)
        \\
        & \geq \liminf_{k \to \infty} \left( \bigl| {\ts \alpha(\bar{\eta}_{n_k})} \gamma_{\bar{\varepsilon}_{n_k}}(\nabla \check{\theta}_{n_k}) \bigr|_{L^1(\Omega)} +\frac{\bar{\nu}_{n_k}^2}{2} \bigl| \nabla \check{\theta}_{n_k} \bigr|_{[X]^N}^2 +\frac{M_0}{2} \bigl| \check{\theta}_{n_k} \bigr|_{X}^2 \right)
        \\
        & \geq \bigl| {\ts \alpha(\bar{\eta})} \gamma_{\bar{\varepsilon}}(\nabla \check{\theta}) \bigr|_{L^1(\Omega)} +\frac{\bar{\nu}^2}{2} \bigl| \nabla \check{\theta} \bigr|_{[X]^N}^2 +\frac{M_0}{2} \bigl| \check{\theta} \bigr|_{X}^2
        \\
        & = \Phi_{\bar{\nu}, \bar{\varepsilon}}(\bar{\eta}, \check{\theta}).
    \end{align*}

    Meanwhile, if we consider a function $ \hat{\theta} \in D(\Phi_{\bar{\nu}, \bar{\varepsilon}}(\bar{\eta},\cdot{})) = Y $ as in the condition of optimality, then assumptions \eqref{bar-nu-eps}, \eqref{Gconv01}, and (\hyperlink{A6l}{A6}) enable us to choose $ \{ \hat{\theta}_n \}_{n = 1}^\infty :=  \{ \hat{\theta} \} \subset Y $ as the sustaining sequence for the condition of optimality. More precisely,
    \begin{align*}
        \Phi_{\bar{\nu}_n, \bar{\varepsilon}_n}(\bar{\eta}_n, \hat{\theta}) & = \int_\Omega \alpha(\bar{\eta}_n) \gamma_{\bar{\varepsilon}_n}(\nabla \hat{\theta}) \, dx +\frac{\bar{\nu}_n^2}{2} \int_\Omega |\nabla \hat{\theta}|^2 \, dx +\frac{M_0}{2} \int_\Omega |\hat{\theta}|^2 \,dx
        \\
        & \to \int_\Omega \alpha(\bar{\eta}) \gamma_{\bar{\varepsilon}}(\nabla \hat{\theta}) \, dx +\frac{\bar{\nu}^2}{2} \int_\Omega |\nabla \hat{\theta}|^2 \, dx +\frac{M_0}{2} \int_\Omega |\hat{\theta}|^2 \,dx
        \\
        & = \Phi_{\bar{\nu}, \bar{\varepsilon}}(\bar{\eta}, \hat{\theta}) \mbox{ as $ n \to \infty $.}
    \end{align*}

    Second, we assume the convergence of optimality \eqref{Gconv02} to verify the strong convergence of the sequence $ \{ \bar{\theta}_n \}_{n = 1}^\infty $ in $ Y $. Given \eqref{bar-nu-eps}, \eqref{Gconv01}, (\hyperlink{A6l}{A6}), and Lions' lemma (cf. \cite[Lemma 1.3 on page 12]{MR0259693}), we observe that
    \begin{align}\label{KenM04}
            \bar{\theta}_n \to & \bar{\theta} \mbox{ in $ X $, weakly in $ Y $, and in the pointwise sense, a.e. in $ \Omega $,}
        \\
        & \hspace{11ex}\mbox{as $ n \to \infty $, by considering a subsequence.}
        \nonumber
    \end{align}
Bearing in mind \eqref{KenM04} and applying a computation similar to \eqref{*ken}, we deduce the following: 
\begin{align}\label{KenM02}
 \liminf_{n \to \infty} \bigl| \alpha (\bar{\eta}_n) \gamma_{\bar{\varepsilon}_n}(\nabla \bar{\theta}_n)\bigr|_{L^1(\Omega)} \geq \bigl| \alpha (\bar{\eta}) \gamma_{\bar{\varepsilon}}(\nabla \bar{\theta})\bigr|_{L^1(\Omega)}.  
\end{align}
Additionally, as a direct consequence of \eqref{KenM04}, it follows that
\begin{align}\label{KenM02-01}
\ds \liminf_{n \to \infty} \frac{\bar{\nu}_n^2}{2} \bigl| \nabla \bar{\theta}_n \bigr|_{[X]^N}^2 \geq \frac{\bar{\nu}^2}{2} \bigl| \nabla \bar{\theta} \bigr|_{[X]^N}^2 \mbox{ and } \lim_{n \to \infty} \frac{M_0}{2} \bigl| \bar{\theta}_n \bigr|_{X}^2 = \frac{M_0}{2} \bigl| \bar{\theta} \bigr|_{X}^2.
\end{align}
Because of \eqref{Gconv02}, \eqref{KenM02}, and \eqref{KenM02-01}, we obtain the following: 
    \begin{align}\label{KenM03}
        & \begin{cases}
            \ds \lim_{n \to \infty} \bigl| {\ts \alpha(\bar{\eta}_n}) \gamma_{\bar{\varepsilon}_n}(\nabla \bar{\theta}_n) \bigr|_{L^1(\Omega)} = \bigl| {\ts \alpha(\bar{\eta}}) \gamma_{\bar{\varepsilon}}(\nabla \bar{\theta}) \bigr|_{L^1(\Omega)},
            \\[2ex]
            \ds \lim_{n \to \infty} \frac{\bar{\nu}_n^2}{2} \bigl| \nabla \bar{\theta}_n \bigr|_{[X]^N}^2 = \frac{\bar{\nu}^2}{2} \bigl| \nabla \bar{\theta} \bigr|_{[X]^N}^2.
        \end{cases}
    \end{align}
    We verify the strong convergence of $ \{ \bar{\theta}_n \}_{n = 1}^\infty $ in $ Y $ as a consequence of  the convergences as in \eqref{bar-nu-eps}, \eqref{KenM04}--\eqref{KenM03}, and the uniform convexity of the $L^2$-based norm. 

    Thus, we complete the proof of this lemma. 
        \qed
}

    \begin{lem}\label{Lem3-ContDep}
        Under assumptions (\hyperlink{A1l}{A1})--(\hyperlink{A7l}{A7}), let $ \bar{\nu} \in (0, \nu_0 + 1] $, $ \{ \bar{\nu}_n \}_{n = 1}^\infty \subset (0, \nu_0 + 1] $, $ \bar{\varepsilon} \in (0, 1) $, and $ \{\bar{\varepsilon}_n\}_{n = 1}^\infty \subset (0, 1) $ be as in \eqref{bar-nu-eps}. Let $ [\eta_0, \theta_0] \in  [Y]^2  $ be an initial pair such that
        \begin{align*}
            & [\eta_0, \theta_0] \in [L^\infty(\Omega)]^2  \mbox{ and } |\eta_0|_{L^\infty(\Omega)} \vee |\theta_0|_{L^\infty(\Omega)} \leq R_0,
        \end{align*}
        and let $ \bigl[ \{\eta_i\}_{i = 1}^m, \{\theta_i\}_{i = 1}^m \bigr] $ be the solution to (S)$_{\bar{\nu}, \bar{\varepsilon}}^{m}$ for the initial pair $ [\eta_0, \theta_0] \in [Y]^2 $ and forcing pair $ \bigl[ \{ u_i^{m} \}_{i = 1}^m, \{ v_i^{m} \}_{i = 1}^m \bigr] \in [L^\infty(\Omega)]^{m \times 2} $ as in Remark \ref{uv.P}. 
        Additionally, let $ \bigl\{ [\eta_{n, 0}, \theta_{n, 0}] \bigr\}_{n = 1}^\infty \subset [Y]^2 $ be a sequence of initial pairs and let $ \bigl\{ \bigl[ \{u_{n, i}^{m}\}_{i = 1}^m, \{v_{n, i}^{m}\}_{i = 1}^m \bigr] \bigr\}_{n = 1}^\infty \subset [L^\infty(\Omega)]^{m \times 2} $ be a sequence of forcing pairs such that
\begin{subequations}
       \begin{align}\label{ShoM00-01}
            & \begin{cases}
                \bigl\{ \bigl[ \{u_{n, i}^{m}\}_{i = 1}^m, \{v_{n, i}^{m}\}_{i = 1}^m \bigr] \bigr\}_{n = 1}^\infty \subset Z(R_0), ~ \bigl\{ [\eta_{n, 0}, \theta_{n, 0}] \bigr\}_{n = 1}^\infty \subset [L^\infty(\Omega)]^2, 
                \\[1ex]
                \ds\sup_{n \in \N} \bigl\{ |\eta_{n, 0}|_{L^\infty(\Omega)} \vee |\theta_{n, 0}|_{L^\infty(\Omega)} \bigr\} \leq R_0,
            \end{cases}
        \end{align}
        and
        \begin{align}\label{ShoM00-02}
 & \begin{cases}
            [\{ u_{n, i}^{m}\}_{i=1}^m, \{ v_{n, i}^{m} \}_{i=1}^m] \to [\{u_i^{m} \}_{i=1}^m, \{ v_i^{m}\}_{i=1}^m] \mbox{ in } [X]^{m \times 2}, \\[1ex]
            [\eta_{n, 0}, \theta_{n, 0}] \to [ \eta_0, \theta_0] \mbox{ weakly in $ [Y]^2 $}
\end{cases}
 \mbox{  as $ n \to \infty $.}
\end{align}
\end{subequations}
        Additionally, for any $ n \in \N $, let $ \bigl[ \{\eta_{n, i}\}_{i = 1}^m, \{\theta_{n, i}\}_{i = 1}^m \bigr] $ be the solution to (S)$_{\bar{\nu}_n, \bar{\varepsilon}_n}^{m}$ for the initial pair $ [\eta_{n, 0}, \theta_{n, 0}] \in [Y]^2 $ and forcing pair $ \bigl[ \{ u_{n, i}^{m} \}_{i = 1}^m, \{ v_{n, i}^{m} \}_{i = 1}^m \bigr] \in [L^\infty(\Omega)]^{m \times 2} $. Then, it holds that
        \begin{align}\label{Sho00}
            & \bigl[ \{ \eta_{n, i} \}_{i = 1}^m, \{ \theta_{n, i}\}_{i = 1}^m \bigr] \to \bigl[ \{ \eta_i \}_{i = 1}^m, \{ \theta_i \}_{i = 1}^m \bigr] \mbox{ in $  [Y]^{m \times 2}  $,}
            \\[1ex]
            & \hspace{7ex} \mbox{and weakly-$*$ in $ [L^\infty(\Omega)]^{m \times 2} $, as $ n \to \infty $.}
            \nonumber
        \end{align}
    \end{lem}

\paragraph{Proof of Lemma \ref{Lem3-ContDep}.}{
    First, we note that assumption \eqref{ShoM00-01} and Lemma \ref{Lem3-01} lead to
    \begin{align}\label{kenken01}
        & \ds\sup_{n \in \N} \max_{1 \leq i \leq m} \left\{ |\eta_{n, i}^m|_{L^\infty(\Omega)} \vee |\theta_{n, i}^m|_{L^\infty(\Omega)} \right\} \leq R_0.
    \end{align}
    Additionally, as a consequence of \eqref{ShoM00-02} and Lemma \ref{Lem-Gamma01}, it follows that
    \begin{align}\label{kenkenMosco}
        & \Phi_{\bar{\nu}_n, \bar{\varepsilon}_n} \to \Phi_{\bar{\nu}, \bar{\varepsilon}} ~\mbox{on $ X $, in the sense of Mosco, as $ n \to \infty $.}
    \end{align}

Next, we fix $n \in \mathbb{N}$ and sum the energy estimates in \eqref{energyEst01}. Then, 
\begin{align}\label{ShoM01}
\ds \frac{1}{2\tau}\sum_{j=1}^i  \Bigl(|\eta_{n, j} - \eta_{n, j-1}|_X^2  & +  |\sqrt{\alpha_0(\eta_{n, j-1})}(\theta_{n, j} - \theta_{n, j-1})|_X^2 \Bigr)  \nonumber \\
    + \mathcal{F}_{\bar{\nu}_n, \bar{\varepsilon}_n}  (\eta_{n, i}, \theta_{n, i}) & \leq  \mathcal{F}_{\bar{\nu}_n, \bar{\varepsilon}_n}(\eta_{n, 0}, \theta_{n, 0}) + i\tau  R_0^2 \mathcal{L}^N(\Omega) \left( 1 +\frac{1}{2 \delta_*^2} \right)  \nonumber \\
    & \leq  \mathcal{F}_{\bar{\nu}_n, \bar{\varepsilon}_n}(\eta_{n, 0}, \theta_{n, 0}) +T  R_0^2 \mathcal{L}^N(\Omega) \left( 1 +\frac{1}{2 \delta_*^2} \right) \\
& \quad  \mbox{ for any } i = 1, 2, 3, \ldots, m. \nonumber
\end{align} 
Bearing in mind \eqref{phi.nu-eps}, \eqref{free.en}, \eqref{bar-nu-eps}, and \eqref{ShoM00-02}, we determine the positive constant $C_8$ independent of $m $ and $ n$ such that
\begin{align}\label{ShoM02}
\ds  \mathcal{F}_{\bar{\nu}_n, \bar{\varepsilon}_n}(\eta_{n, 0}, \theta_{n, 0}) & = \frac{\kappa^2}{2}|\nabla \eta_{n, 0}|_{[X]^N}^2 + \int_\Omega G(\eta_{n, 0})\, dx \nonumber \\
                                & \, \quad + \int_\Omega \alpha (\eta_{n, 0}) \gamma_{\bar{\varepsilon}_n}(\nabla \theta_{n, 0})\, dx + \frac{\bar{\nu}_n^2}{2}|\nabla \theta_{n, 0}|_{[X]^N}^2 + \frac{M_0}{2}|\theta_{n, 0}|_X^2 \nonumber \\
                                & \leq  \frac{\kappa^2}{2}|\eta_{n, 0}|_Y^2 + |G|_{C([-R_0, R_0])} \mathcal{L}^N(\Omega) \nonumber \\
 & \, \quad + \mathcal{L}^N(\Omega)|\alpha|_{L^\infty(-R_0, R_0)}\bar{\varepsilon}_n + \sqrt{\mathcal{L}^N(\Omega)}|\alpha|_{L^\infty(-R_0, R_0)}|\nabla \theta_{n, 0}|_{[X]^N} \nonumber \\
& \, \quad + \frac{(1 + \nu_0)^2 \vee M_0}{2}|\theta_{n, 0}|_Y^2 \nonumber \\
& \leq C_8. 
\end{align}
Taking into account (\hyperlink{A3l}{A3}), \eqref{ShoM01}, \eqref{ShoM02}, and the estimate
\begin{align*}
    & \hspace{15ex}\frac{1}{4T} (|\eta_{n, i}|_X^2 + \delta_*|\theta_{n, i}|_X^2) -\frac{1}{2T}(|\eta_{n, 0}|_X^2 +  \delta_*|\theta_{n, 0}|_X^2)  
    \\
    & \leq \frac{1}{2\tau}\sum_{j=1}^i  \left(|\eta_{n, j} - \eta_{n, j-1}|_X^2 +  \bigl|{\textstyle\sqrt{\alpha_0(\eta_{n, j-1})}}(\theta_{n, j} - \theta_{n, j-1}) \bigr|_X^2 \right), ~ i = 1, 2, 3, \dots, m,
\end{align*}
we compute the following: 
\begin{align}\label{ShoM03}
 \ds & \frac{1}{4T}(|\eta_{n, i}|_X^2 + \delta_*|\theta_{n, i}|_X^2) + \mathcal{F}_{\bar{\nu}_n, \bar{\varepsilon}_n}  (\eta_{n, i}, \theta_{n, i}) \nonumber 
    \\
  & \leq \frac{1}{2T}(|\eta_{n, 0}|_X^2 +  \delta_*|\theta_{n, 0}|_X^2) \nonumber 
    \\
& \quad + \frac{1}{2\tau}\sum_{j=1}^i  \left(|\eta_{n, j} - \eta_{n, j-1}|_X^2 +  |\sqrt{\alpha_0(\eta_{n, j-1})}(\theta_{n, j} - \theta_{n, j-1})|_X^2 \right) + \mathcal{F}_{\bar{\nu}_n, \bar{\varepsilon}_n}  (\eta_{n, i}, \theta_{n, i}) \nonumber 
    \\
 & \leq \frac{1}{2T}(|\eta_{n, 0}|_X^2 + \delta_*|\theta_{n, 0}|_X^2) + \mathcal{F}_{\bar{\nu}_n, \bar{\varepsilon}_n}(\eta_{n, 0}, \theta_{n, 0}) +T R_0^2 \mathcal{L}^N(\Omega) \left( 1 +\frac{1}{2 \delta_*^2} \right) \nonumber 
    \\
    & \leq \frac{1}{2T}\sup_{n \in \N}(|\eta_{n, 0}|_X^2 + \delta_*|\theta_{n, 0}|_X^2) + C_6 + T  R_0^2 \mathcal{L}^N(\Omega) \left( 1 +\frac{1}{2 \delta_*^2} \right) \nonumber 
    \\
 & =: C_9 \mbox{ for any } i = 1, 2, 3, \ldots, m.
\end{align}
Note that $C_9$ can be assumed a positive constant independent of $m$ and $ n$. 

As a consequence of \eqref{free.en}, \eqref{ShoM01}, and \eqref{ShoM03}, we observe that 
\begin{align}\label{ShoM04}
    \{[\eta_n, \theta_n ]\}_{n = 1}^\infty = \{[\{\eta_{n, i}\}_{i=1}^m, \{\theta_{n, i}\}_{i=1}^m] \}_{n = 1}^\infty \mbox{ is bounded in } [Y]^{ m \times 2} \cap [L^\infty(\Omega)]^{m \times 2}. 
\end{align}
Because of \eqref{ShoM04} and the compact embedding $Y \subset X$, there exist $\{ n_k \}_{k = 1}^\infty \subset \{n\}$ and $[\bar{\eta}, \bar{\theta}] = [\{\bar{\eta}_i \}_{i=1}^m, \{\bar{\theta}_i \}_{i=1}^m] \in [Y]^{ m \times 2 }$ such that 
\begin{align}\label{ShoM05}
    [ \eta_{n_k}, \theta_{n_k}] \to & [\bar{\eta}, \bar{\theta}] \mbox{ in }  [X]^{m \times 2},  \mbox{ weakly in }  [Y]^{m \times 2}, \mbox{ weakly-$*$ in $ [L^\infty(\Omega)]^{m \times 2} $,} \nonumber \\
& \mbox{ and in the pointwise sense, a.e. in } \Omega, \mbox{ as } k \to \infty. 
\end{align}
 
Next, we verify that the limit $[\bar{\eta}, \bar{\theta}]$ is the solution to \hyperlink{(S)$_{\nu, \varepsilon}^{m}$}{(S)$_{\nu, \varepsilon}^{m}$}.  
We fix $k \in \mathbb{N}$. 
Then, the solution $[\eta_{n_k}, \theta_{n_k} ] $ admits \eqref{3tau} as follows:
\begin{equation}\label{ShoM06}
\left\{ \parbox{13.5cm}{
    $\displaystyle \frac{1}{\tau} (\eta_{n_k, i} - \eta_{n_k, i-1}, \varphi)_X + \kappa^2( \nabla \eta_{n_k, i}, \nabla \varphi)_{[X]^N} + \int_{\Omega} g(\eta_{n_k, i})\varphi\, dx 
\\[1.5ex]
 + \int_\Omega \alpha'(\eta_{n_k, i}) \gamma_{\bar{\varepsilon}_{n_k}}(\nabla \theta_{n_k, i})\varphi\, dx = (u_{n_k, i}^{m}, \varphi)_X $ for any $\varphi \in Y$
\\[1ex]
$\mbox{ for any } i = 1, 2, 3, \ldots , m;$
}\right. 
\end{equation}
\begin{equation}\label{ShoM07}
\left\{\parbox{13.5cm}{
    $ \displaystyle \frac{1}{\tau} (\alpha_{0}(\eta_{n_k, i-1})( \theta_{n_k, i} - \theta_{n_k, i-1}), \psi)_X + \bar{\nu}_{n_k}^2(\nabla\theta_{n_k, i}, \nabla \psi)_{[X]^N} 
\\[1.5ex] 
+ \int_\Omega \alpha(\eta_{n_k, i-1})\nabla \gamma_{\bar{\varepsilon}_{n_k}} (\nabla \theta_{n_k, i}) \cdot \nabla \psi \, dx  + M_0(\theta_{n_k, i}, \psi)_X = (v_{n_k, i}^{m}, \psi)_X $\\[1ex]
 for any $\psi \in Y$ $ \mbox{ for any } i = 1, 2, 3, \ldots , m$.
}\right. 
\end{equation}

We fix $i \in \{ 1, 2, 3, \ldots, m \}$ and apply Remark \ref{Shorem01} to \eqref{ShoM07}: 
\begin{align*}
 \displaystyle  \Bigl( \frac{1}{\tau} (\alpha_{0}(\eta_{n_k, i -1})& ( \theta_{n_k, i}  - \theta_{n_k, i-1})) - v_{n_k, i}^{m} , \theta_{n_k, i} - \psi \Bigr)_X \\
 & + \Phi_{\bar{\nu}_{n_k}, \bar{\varepsilon}_{n_k}}(\eta_{n_k, i-1}, \theta_{n_k, i}) - \Phi_{\bar{\nu}_{n_k}, \bar{\varepsilon}_{n_k}}(\eta_{n_k, {i-1}}, \psi) \leq 0;
\end{align*}
that is,
\begin{align}\label{ShoM08}
\left[\theta_{n_k, i},  v_{n_k, i}^{m} - \frac{1}{\tau} (\alpha_{0}(\eta_{n_k, i-1})( \theta_{n_k, i} - \theta_{n_k, i-1}) \right] \in \partial \Phi_{\bar{\nu}_{n_k}, \bar{\varepsilon}_{n_k}}(\eta_{n_k, i-1}, \cdot) \mbox{ in } X \times X. 
\end{align}
Bearing in mind (\hyperlink{A3l}{A3}), \eqref{ShoM05}, and \eqref{ShoM08}, and applying Lemma \ref{Lem-Gamma01}  and (\hyperlink{Fact1}{Fact\,1}), we obtain
\begin{align*}
\left[\bar{\theta}_i,  v_i^{m} - \frac{1}{\tau} (\alpha_{0}(\bar{\eta}_{i-1})( \bar{\theta}_i - \bar{\theta}_{i-1}) \right] \in \partial \Phi_{\bar{\nu}, \bar{\varepsilon}}(\bar{\eta}_{i-1}, \cdot) \mbox{ in } X \times X; 
\end{align*}
that is,
\begin{align*}
\displaystyle & \Phi_{\bar{\nu}, \bar{\varepsilon}}(\bar{\eta}_{i-1}, \bar{\theta}_i) - \Phi_{\bar{\nu}, \bar{\varepsilon}} (\bar{\eta}_{i-1}, \psi) \leq  \left( v_i^{m} - \frac{1}{\tau} (\alpha_{0}(\bar{\eta}_{i-1})( \bar{\theta}_i - \bar{\theta}_{i-1}), \bar{\theta}_i - \psi \right)_X.
\end{align*}
This inequality and Remark \ref{Shorem01} lead to 
\begin{equation}\label{ShoM09}
\left\{\parbox{13.5cm}{
    $ \displaystyle \frac{1}{\tau} (\alpha_{0}(\bar{\eta}_{i-1})(  \bar{\theta}_i -  \bar{\theta}_{i-1}), \psi)_X + \bar{\nu}^2(\nabla \bar{\theta}_i, \nabla \psi)_{[X]^N} 
\\[1.5ex] 
+ \int_\Omega \alpha(\bar{\eta}_{i-1})\nabla \gamma_{\bar{\varepsilon}} (\nabla  \bar{\theta}_i) \cdot \nabla \psi \, dx  + M_0( \bar{\theta}_i, \psi)_X = (v_i^{m}, \psi)_X $ for any $\psi \in Y$
\\[1ex]
$ \mbox{ for any } i = 1, 2, 3, \ldots , m$.
}\right. 
\end{equation}

Next, we fix $i = 1, 2, 3, \ldots, m$. 
Then, applying Remark \ref{Shorem01} to \eqref{ShoM07} and putting $\psi = \bar{\theta}_i $ in Remark \ref{Shorem01}, we obtain 
\begin{align*}
 \displaystyle  \Bigl( \frac{1}{\tau} (\alpha_{0}(\eta_{n_k, i -1})& ( \theta_{n_k, i}  - \theta_{n_k, i-1})) - v_i^{m} , \theta_{n_k, i} - \bar{\theta}_i \Bigr)_X \\
 & + \Phi_{\bar{\nu}_{n_k}, \bar{\varepsilon}_{n_k}}(\eta_{n_k, i-1}, \theta_{n_k, i}) - \Phi_{\bar{\nu}_{n_k}, \bar{\varepsilon}_{n_k}}(\eta_{n_k, i-1}, \bar{\theta}_i) \leq 0;
\end{align*}
that is, 
\begin{align}\label{ShoM10}
 \displaystyle   \Bigl( \frac{1}{\tau} (\alpha_{0}(\eta_{n_k, i -1})& ( \theta_{n_k, i}  - \theta_{n_k, i-1})) - v_i^{m} , \theta_{n_k, i} - \bar{\theta}_i \Bigr)_X \nonumber \\
 & + \Phi_{\bar{\nu}_{n_k}, \bar{\varepsilon}_{n_k}}(\eta_{n_k, i-1}, \theta_{n_k, i}) \leq \Phi_{\bar{\nu}_{n_k}, \bar{\varepsilon}_{n_k}}(\eta_{n_k, i-1}, \bar{\theta}_i). 
\end{align}
Because of (\hyperlink{A3l}{A3}), \eqref{ShoM05}, and \eqref{ShoM10}, we observe that
\begin{align}\label{ShoM11}
 \limsup_{k \to \infty} \Phi_{\bar{\nu}_{n_k}, \bar{\varepsilon}_{n_k}}(\eta_{n_k, i-1}, \theta_{n_k, i}) \leq \limsup_{k \to \infty}\Phi_{\bar{\nu}_{n_k}, \bar{\varepsilon}_{n_k}}(\eta_{n_k, i-1}, \bar{\theta}_i).
\end{align}
The convergences \eqref{bar-nu-eps} and \eqref{ShoM05} lead to 
\begin{align*}
\ds \limsup_{k \to \infty} & \int_\Omega \alpha(\eta_{n_k, i-1}) \gamma_{\bar{\varepsilon}_{n_k}}(\nabla \bar{\theta}_i)\, dx \\
& \leq \int_\Omega \alpha(\bar{\eta}_{i-1}) \gamma_{\bar{\varepsilon}}(\nabla \bar{\theta}_i)\, dx + \limsup_{k \to \infty} \Bigl( \int_\Omega |\alpha(\eta_{n_k, i-1}) - \alpha(\bar{\eta}_{i-1})| \gamma_{\bar{\varepsilon}}(\nabla \bar{\theta}_i)\, dx  \\
& \, \quad + \int_\Omega \alpha(\eta_{n_k, i-1})( \gamma_{\bar{\varepsilon}_{n_k}}(\nabla \bar{\theta}_i) - \gamma_{\bar{\varepsilon}}(\nabla \bar{\theta}_i)) \, dx \Bigr)  \\
& \leq \int_\Omega \alpha(\bar{\eta}_{i-1}) \gamma_{\bar{\varepsilon}}(\nabla \bar{\theta}_i)\, dx + |\alpha'|_{L^\infty(-R_0, R_0)} \left( \lim_{k \to \infty}|\eta_{n_k, i-1} - \bar{\eta}_{i-1}|_X \right) \bigl| \bar{\varepsilon} + |\nabla \bar{\theta}_i|\bigr|_X \\
& \, \quad + \sup_{k \in \mathbb{N}}|\alpha (\eta_{n_k, i-1})|_{L^\infty(\Omega)} \lim_{k \to \infty}|\bar{\varepsilon}_{n_k} - \bar{\varepsilon}|_{L^1(\Omega)}  \\
& = \int_\Omega \alpha(\bar{\eta}_{i-1}) \gamma_{\bar{\varepsilon}}(\nabla \bar{\theta}_i)\, dx.
\end{align*}
Hence,
\begin{align}\label{ShoM12}
\ds \limsup_{k \to \infty}& \Phi_{\bar{\nu}_{n_k}, \bar{\varepsilon}_{n_k}}(\eta_{n_k, i-1}, \bar{\theta}_i) \nonumber \\
 & \leq  \lim_{k \to \infty}\frac{\bar{\nu}_{n_k}^2}{2}|\nabla \bar{\theta}_i|_{[X]^N}^2 + \limsup_{k \to \infty} \int_\Omega \alpha(\eta_{n_k, i-1}) \gamma_{\bar{\varepsilon}_{n_k}}(\nabla \bar{\theta}_i)\, dx + \frac{M_0}{2}|\bar{\theta}_i|_{X}^2 \nonumber \\
& \leq \frac{\bar{\nu}^2}{2}|\nabla \bar{\theta}_i|_{[X]^N}^2 + \int_\Omega \alpha(\bar{\eta}_{i-1}) \gamma_{\bar{\varepsilon}}(\nabla \bar{\theta}_i)\, dx + \frac{M_0}{2}|\bar{\theta}_i|_{X}^2 \nonumber \\
& = \Phi_{\bar{\nu}, \bar{\varepsilon}}(\bar{\eta}_{i-1}, \bar{\theta}_i). 
\end{align}
Taking into account \eqref{kenkenMosco}, \eqref{ShoM11}, and \eqref{ShoM12}, \begin{align*}
 \ds \Phi_{\bar{\nu}, \bar{\varepsilon}}(\bar{\eta}_{i-1}, \bar{\theta}_i) & \leq \liminf_{k \to \infty} \Phi_{\bar{\nu}_{n_k}, \bar{\varepsilon}_{n_k}}(\eta_{n_k, i-1}, \theta_{n_k, i}) \leq \limsup_{k \to \infty} \Phi_{\bar{\nu}_{n_k}, \bar{\varepsilon}_{n_k}}(\eta_{n_k, i-1}, \theta_{n_k, i})\\
&  \leq   \limsup_{k \to \infty} \Phi_{\bar{\nu}_{n_k}, \bar{\varepsilon}_{n_k}}(\eta_{n_k, i-1}, \bar{\theta}_i) \leq \Phi_{\bar{\nu}, \bar{\varepsilon}}(\bar{\eta}_{i-1}, \bar{\theta}_i).
\end{align*}
Hence,
\begin{align}\label{ShoM13}
 \ds \lim_{k \to \infty} \Phi_{\bar{\nu}_{n_k}, \bar{\varepsilon}_{n_k}}(\eta_{n_k, i-1}, \theta_{n_k, i})  =  \Phi_{\bar{\nu}, \bar{\varepsilon}}(\bar{\eta}_{i-1}, \bar{\theta}_i). 
\end{align}
From the convergence \eqref{ShoM05}, \eqref{ShoM13}, and Lemma \ref{Lem-Gamma01}, we observe that 
\begin{align}\label{ShoM14}
 \theta_{n_k} \to \bar{\theta} \mbox{ in } [Y]^m \mbox{ as } k \to \infty. 
\end{align}

Now, we fix arbitrary $i \in \{ 1, 2, 3, \ldots, n \} $. 
Note that it follows from \eqref{ShoM14} that
\begin{align}\label{ShoM15}
 |\gamma_{\bar{\varepsilon}_{n_k}} (\nabla \theta_{n_k, i})&  - \gamma_{\bar{\varepsilon}} (\nabla \bar{\theta}_i)|_X \leq |\gamma_{\bar{\varepsilon}_{n_k}} (\nabla \theta_{n_k, i}) - \gamma_{\bar{\varepsilon}_{n_k}} (\nabla \bar{\theta}_i)|_X + |\gamma_{\bar{\varepsilon}_{n_k}} (\nabla \bar{\theta}_i) - \gamma_{\bar{\varepsilon}} (\nabla \bar{\theta}_i)|_X \nonumber \\
& \leq |\theta_{n_k, i} - \bar{\theta}_i|_Y + |\bar{\varepsilon}_{n_k} - \bar{\varepsilon}|_{L^1(\Omega)} \to 0 \mbox{ as } k \to \infty. 
\end{align}
The convergence \eqref{ShoM05} and the assumptions lead to
\begin{align}\label{ShoM16}
            \ds \alpha'( \eta_{n_k, i})\varphi \to  \alpha'(\bar{\eta}_i) \varphi, \mbox{ in the pointwise sense, a.e. in $ \Omega $, as $ k \to \infty $} 
\end{align}
as a result of considering a subsequence, if necessary. 
As a consequence of (\hyperlink{A3l}{A3}), \eqref{ShoM05}, \eqref{ShoM16}, and Lebesgue's dominated convergence theorem, we observe that
\begin{align}\label{ShoM17}
 \alpha'(\eta_{n_k, i})\varphi \to  \alpha'(\bar{\eta}_i)\varphi \mbox{ in } X \mbox{ as } k \to \infty.
\end{align}
From (\hyperlink{A4l}{A4}), \eqref{ShoM05}, \eqref{ShoM15}, and \eqref{ShoM17}, it follows that
\begin{subequations}\label{ShoM18}
\begin{align}\label{ShoM18-1}\noeqref{ShoM18-1}
 \displaystyle \frac{1}{\tau} (\eta_{n_k, i} - \eta_{n_k, i-1}, \varphi)_X \to \frac{1}{\tau} (\bar{\eta}_i - \bar{\eta}_{i-1}, \varphi)_X \mbox{ as } k \to \infty,
\end{align}
\begin{align}\label{ShoM18-2}\noeqref{ShoM18-2}
\kappa^2( \nabla \eta_{n_k, _i}, \nabla \varphi)_{[X]^N} \to \kappa^2( \nabla \bar{\eta}_i, \nabla \varphi)_{[X]^N} \mbox{ as } k \to \infty,
\end{align}
\begin{align}\label{ShoM18-3}\noeqref{ShoM18-3}
g(\eta_{n_k, i}) \to  g(\bar{\eta}_i) \mbox{ in } X \mbox{ as } k \to \infty,
\end{align}
\begin{align}\label{ShoM18-4}\noeqref{ShoM18-4}
\ds \int_\Omega \alpha'(\eta_{n_k, i}) \gamma_{\bar{\varepsilon}_{n_k}}(\nabla \theta_{n_k, i})\varphi\, dx &\, = ( \gamma_{\bar{\varepsilon}_{n_k}}(\nabla \theta_{n_k, i}), \alpha'(\eta_{n_k, i}) \varphi )_X \nonumber \\ 
& \, \to ( \gamma_{\bar{\varepsilon}}(\nabla \bar{\theta}_{i}), \alpha'(\bar{\eta}_i) \varphi )_X \nonumber \\
& \, = \int_\Omega \alpha'(\bar{\eta}_i) \gamma_{\bar{\varepsilon}}(\nabla \bar{\theta}_{i})\varphi\, dx \mbox{ as } k \to \infty.
\end{align}
\end{subequations}
The convergences \eqref{ShoM05} and \eqref{ShoM18}, and \eqref{ShoM06} lead to 
\begin{equation}\label{ShoM19}
\left\{ \parbox{13.5cm}{
    $\displaystyle \frac{1}{\tau} (\bar{\eta}_i - \bar{\eta}_{i-1}, \varphi)_X + \kappa^2( \nabla \bar{\eta}_i, \nabla \varphi)_{[X]^N} + \int_{\Omega} g(\bar{\eta}_i)\varphi\, dx 
\\[1.5ex]
 + \int_\Omega \alpha'(\bar{\eta}_i) \gamma_{\bar{\varepsilon}}(\nabla \bar{\theta}_{i})\varphi\, dx = (u_i^{m}, \varphi)_X $ for any $\varphi \in Y$
\\[1ex]
$\mbox{ for any } i = 1, 2, 3, \ldots , m.$
}\right. 
\end{equation}
\eqref{ShoM09} and \eqref{ShoM19} imply that the limit $ [\bar{\eta}, \bar{\theta}] $ coincides with the unique solution $ [\eta, \theta] $ to the system (S)$_{\bar{\nu}, \bar{\varepsilon}}^{m}$. Therefore, we can say that
\begin{itemize}
    \item[$(**)$]the convergences \eqref{ShoM05} and \eqref{ShoM14} are true, even if $ \{n_k\} $ and $ [\bar{\eta}, \bar{\theta}] $ are replaced by $ \{n\} $ and $ [\eta, \theta] $, respectively.
\end{itemize}

Finally, we verify the strong convergence of the sequence $\{ \eta_{n} \}_{n = 1}^\infty$ to $\eta$ in $Y$. 
As a direct consequence of \eqref{ShoM05} and $(**)$, we observe that 
\begin{align}\label{ShoM20}
 \liminf_{n \to \infty} \kappa^2|\nabla \eta_{n, i}|_{[X]^N}^2 \geq \kappa^2|\nabla {\eta}_i|_{[X]^N}^2.
\end{align}
Bearing \eqref{ShoM05}, \eqref{ShoM06}, and \eqref{ShoM17}--\eqref{ShoM19} in mind, we compute the following: 
\begin{align}\label{ShoM21}
 \limsup_{n \to \infty} \kappa^2|\nabla \eta_{n, i}|_{[X]^N}^2 & = \limsup_{k \to \infty} \Bigl( (u_{n, i}^{m}, \eta_{n, i})_X - \frac{1}{\tau} (\eta_{n, i} - \eta_{n, i-1}, \eta_{n, i})_X \nonumber 
\\
 & \quad - \int_{\Omega} g(\eta_{n, i})\eta_{n, i}\, dx - \int_\Omega \alpha'(\eta_{n, i}) \gamma_{\bar{\varepsilon}_{n}}(\nabla \theta_{n, i})\eta_{n, i}\, dx\Bigr) \nonumber
\\
& \leq  (u_i^{m}, {\eta}_i)_X - \frac{1}{\tau} ({\eta}_i - {\eta}_{i-1}, {\eta}_i )_X \nonumber 
\\
& \quad - \int_{\Omega} g({\eta}_i){\eta}_i\, dx - \int_\Omega \alpha'({\eta}_i) \gamma_{{\varepsilon}}(\nabla {\theta}_{i}){\eta}_i \, dx \nonumber
\\ 
& = \kappa^2|\nabla {\eta}_i|_{[X]^N}^2. 
\end{align}
The inequalities \eqref{ShoM20} and \eqref{ShoM21} lead to 
\begin{align}\label{ShoM22} 
 \lim_{n \to \infty} \frac{\kappa^2}{2}| \nabla \eta_{n, i} |_{[X]^N}^2 = \frac{\kappa^2}{2}|\nabla {\eta}_i|_{[X]^N}^2.
\end{align} 
We verify the strong convergence of $\{ \eta_{n_k} \}_{k=1}^\infty $ in $Y$ as a consequence of \eqref{ShoM05}, \eqref{ShoM22}, $ (**) $, and the uniform convexity of the $L^2$-based norm. 

Thus, we conclude the proof of Lemma \ref{Lem3-ContDep}. 
\qed
}

\paragraph{Proof of Theorem \ref{Thm.1}.}{
    We fix $ \nu \in (0, 1) $ and define a class $\mathbb{K}_\nu^m \subset [X]^2$ by letting 
\begin{equation}\label{thm1-01}
    \mathbb{K}_\nu^m := \left\{ \begin{array}{l|l}
         [\tilde{\eta}, \tilde{\theta}] \in [Y]^2 & 
         \parbox{7.77cm}{
             $ |\tilde{\eta}|_{L^\infty(\Omega)} \vee |\tilde{\theta}|_{L^\infty(\Omega)} \leq R_0 $, 
             \\[1ex]
             $ \ds |\tilde{\eta}|_X^2 + R_*|\tilde{\theta}|_X^2 +  \kappa^2 \tau|\nabla \tilde{\eta}|_{[X]^N}^2  \leq R_1 $, 
             \\[1ex] 
            $ \kappa^2 |\nabla \tilde{\eta}|_{[X]^N}^2 +\delta_*|D\tilde{\theta}|(\Omega) +\nu^2 |\nabla \tilde{\theta}|_{[X]^N}^2 \leq  R_3  $
    }
    \end{array} \right\},
\end{equation}
where $R_0$ is the positive constant as in \eqref{ken00}, and $ R_1 $ and $ R_3 $ are the positive constants as in Lemma \ref{Lem3-03KS} (a)(b). Then, because of Lemma \ref{Lem3-01}, we can also define an operator $ \mathcal{S} : \mathbb{K}_\nu^m \longrightarrow [X]^2 $, which maps any $ [\eta_0, \theta_0] \in \mathbb{K}_\nu^m $ to the terminal value $ [\eta_m, \theta_m] \in [Y]^2 $ of the solution $ \bigl\{ [\eta_i, \theta_i] \bigr\}_{i = 1}^m \subset [X]^2 $ to the approximating problem (S)$_{\nu, \varepsilon}^{m}$.

The class $ \mathbb{K}_\nu^m $ is a closed convex and bounded set in $ [H^1(\Omega)  \cap L^\infty(\Omega)]^2 $. Hence, we observe that $ \mathbb{K}_\nu^m $ is a compact convex set in $ [X]^2 $. Additionally, by invoking Lemma \ref{Lem3-02} (III) and Lemma \ref{Lem3-03KS}, we observe that $ \mathcal{S}(\mathbb{K}_\nu^m) \subset \mathbb{K}_\nu^m  $; that is, $ \mathcal{S} $ is a selfmapping from $ \mathbb{K}_\nu^m $ to itself. Furthermore, from Lemma \ref{Lem3-ContDep}, we infer the continuity of $ \mathcal{S} $ in the following sense:
\begin{subequations}
\begin{align}
    & \mathcal{S}[\eta_{n, 0}, \theta_{n, 0}] \to \mathcal{S}[\eta_0, \theta_0] \mbox{ in $ [Y]^2 $,}
    \label{contSa}
    \\
    & \mbox{ if } [\eta_0, \theta_0] \in \mathbb{K}_\nu^m, ~ \bigl\{ [\eta_{n, 0}, \theta_{n, 0}] \bigr\}_{n = 1}^\infty \subset \mathbb{K}_\nu^m, \nonumber 
\\
  & \mbox{ and } [\eta_{n, 0}, \theta_{n, 0}] \to [\eta_0, \theta_0] \mbox{ in $ [X]^2 $ as $ n \to \infty $.}
    \label{contSb}
\end{align}
\end{subequations}
In fact, because \eqref{thm1-01} and \eqref{contSb} lead to
\begin{align*}
    & [ \eta_{n, 0}, \theta_{n, 0}] \to [\eta_0, \theta_0] \mbox{ weakly in $ [Y]^2 $ as $ n \to \infty $,}
\end{align*}
the continuous dependence \eqref{contSa} immediately follows as a straightforward consequence of Lemma \ref{Lem3-ContDep}. 

In view of this, we apply Schauder's fixed-point theorem to obtain a fixed point $ [\eta_0^P, \theta_0^P] \in \mathbb{K}_\nu^m$. Indeed, the solution $ \bigl\{ [\eta_i^P, \theta_i^P] \bigr\}_{i = 1}^m \subset [X]^2 $, starting from $ [\eta_0^P, \theta_0^P] \in \mathbb{K}_\nu^m $, is the required periodic solution to the approximating problem (S)$_{\nu, \varepsilon}^{m}$. 

Thus, we conclude the proof of Theorem \ref{Thm.1}.
\qed}

\section{Proof of the Main Theorem}

    Before we prove the Main Theorem, we present the following lemma concerned with the $\Gamma$-convergence of the sequence of governing convex energies on $ L^2(0, T; X) $. 

\begin{lem}\label{Lem3-Gamma}
    Under assumptions (\hyperlink{A1l}{A1})--(\hyperlink{A4l}{A4}) and (\hyperlink{A6l}{A6}),  let $ I \subset (0, T) $ be an open interval, let $ \{ \nu_n \}_{n = 1}^\infty \subset (0, \nu_0 + 1] $, and let $ \{ \varepsilon_n \}_{n = 1}^\infty \subset (0, 1) $ be sequences such that
    \begin{align}\label{bar-nu-eps01} 
        & \nu_n \to \nu_0 ~\mbox{and}~ \varepsilon_n \to 0 ~\mbox{as $ n \to \infty $.}
    \end{align}
    Let $ \widetilde{\eta} \in C(\overline{I}; X) \cap L^\infty(I; Y) \cap L^\infty(Q) $ and $ \{ \widetilde{\eta}_n \}_{n = 1}^\infty \subset L^2(I; Y) $ be such that
    \begin{equation}\label{Gconv00}
    \left\{ \hspace{-3ex} \parbox{9.5cm}{
        \vspace{-2ex}
        \begin{itemize}
            \item$ \{ \widetilde{\eta}_{n}\}_{n = 1}^\infty  $ is bounded in $ L^\infty(I; Y) \cap L^\infty(Q) $, 
            \vspace{-1ex}
            \item$ \widetilde{\eta}_{n} \to \widetilde{\eta} $ in the pointwise sense, a.e. in $ Q $. 
        \end{itemize}
        \vspace{-2ex}
    } \right. 
\end{equation}
    Let $ \widetilde{\Phi}_{0}^I(\widetilde{\eta}, \cdot) : L^2(I; X) \longrightarrow [0, \infty] $ be a proper, l.s.c., and convex function on $ L^2(I; X) $ defined as
    \begin{align*}
        \zeta \in L^2(I; X) ~& \mapsto \widetilde{\Phi}_{0}^I(\widetilde{\eta}, \zeta) := \begin{cases}
            \ds \int_I \Phi_{0}(\widetilde{\eta}(t), \zeta(t)) \, dt, 
            \\
            \qquad \mbox{if $ |D\zeta(\cdot)| \in L^1(I) $, and $ \nu_0 \zeta \in L^2(I; Y) $,}
        \\[1ex]
        \infty, \mbox{ otherwise,}
        \end{cases}
    \end{align*}
    and let $ \{ \widetilde{\Phi}_{\nu_n, \varepsilon_n}^I(\widetilde{\eta}_n, \cdot) \}_{n = 1}^\infty $ be a sequence of proper, l.s.c., and convex functions $ \widetilde{\Phi}_{\nu_n, \varepsilon_n}^I(\widetilde{\eta}_n, \cdot) : L^2(I; X) \longrightarrow [0, \infty] $ defined as
    \begin{align*}
        \zeta \in L^2(I; X) ~& \mapsto \widetilde{\Phi}_{\nu_n, \varepsilon_n}^I(\widetilde{\eta}_n, \zeta) := \begin{cases}
            \ds \int_I \Phi_{\nu_n, \varepsilon_n}(\widetilde{\eta}_n(t), \zeta(t)) \, dt, 
            \\
            \qquad \mbox{if $ \zeta \in L^2(I; Y) $,}
            \\[2ex]
            \infty, \mbox{ otherwise}
        \end{cases}\mbox{for $ n = 1, 2, 3, \dots $.}
    \end{align*}
    Then, the following two items hold:
    \begin{description}
        \item[\textmd{(I)}]$ \widetilde{\Phi}_{\nu_n, \varepsilon_n}^I(\widetilde{\eta}_n, \cdot) \to \widetilde{\Phi}_{0}^I(\widetilde{\eta}, \cdot) $ on $ L^2(I; X) $, in the sense of $ \Gamma $-convergence, as $ n \to \infty $. 
        \item[\textmd{(II)}]For any $ \widetilde{\theta} \in L^2(I; X) $ with $ |D \widetilde{\theta}(\cdot)|(\Omega) \in L^1(I) $, the convergence
        \begin{equation}\label{convOpt}
            \widetilde{\theta}_n \to \widetilde{\theta} \mbox{ in $ C(\overline{I}; X) $, and $ \widetilde{\Phi}_{\nu_n, \varepsilon_n}^I(\widetilde{\eta}_n, \widetilde{\theta}_n) \to {\widetilde{\Phi}_{0}}^I(\widetilde{\eta}, \widetilde{\theta}) $ as $ n \to \infty $}
        \end{equation}
            implies that
            \begin{subequations}\label{kenConv}
           \begin{align}
               &\ds \int_I |D \widetilde{\theta}_n(t)|(\Omega) \, dt \to \int_I |D \widetilde{\theta}(t)|(\Omega) \, dt,
               \label{kenConv01}
               \\[2ex]
               & \hspace{-2ex} \mbox{and }~ \nu_0 \widetilde{\theta}_n \to \nu_0 \widetilde{\theta} \mbox{ in $ L^2(I; Y) $} ~ \mbox{as $ n \to \infty $.}
               \label{kenConv02}
        \end{align}
            \end{subequations}
    \end{description}
\end{lem}

\paragraph{Proof of Lemma \ref{Lem3-Gamma}.}{From assumptions \eqref{Gconv00} and (\hyperlink{A3l}{A3}), we observe that
    \begin{align}\label{gamma00}
        & \begin{cases}
            \alpha(\widetilde{\eta}) \in C(\overline{I}; X) \cap L^\infty(I; Y), ~ \log \alpha(\widetilde{\eta}) \in L^\infty(I \times \Omega),
            \\[1ex]
            \alpha(\widetilde{\eta}_n) \in L^\infty(I; Y), ~ \log \alpha(\widetilde{\eta}_n) \in L^\infty(I \times \Omega) \mbox{ for $ n = 1, 2,3 , \dots $,}
            \\[1ex]
            \alpha(\widetilde{\eta}_n(t)) \to \alpha(\widetilde{\eta}(t)) \mbox{ in $ X $ and weakly in $ Y $ as $ n \to \infty $, for a.e. $ t \in I $.}
        \end{cases}
    \end{align}
    Based on this, we prove item (I). For the condition of the lower bound, we take a function $ \check{\theta} \in L^2(I; X) $ and a sequence $ \{ \check{\theta}_n \}_{n = 1}^\infty \subset L^2(I; X) $ such that
    \begin{align}\label{gamma02}
        \check{\theta}_n \to \check{\theta}; \mbox{ hence $ \nu_n \check{\theta}_n \to \nu_0 \check{\theta} $  in $ L^2(I; X) $ as $ n \to \infty $.}
    \end{align}
    Then, we may suppose the existence of a subsequence $ \{ n_k \}_{k = 1}^\infty \subset \{ n \} $ such that
    \begin{align*}
        \lim_{k \to \infty} \widetilde{\Phi}_{\nu_{n_k}, \varepsilon_{n_k}}^I(\widetilde{\eta}_{n_k}, \check{\theta}_{n_k}) = \liminf_{n \to \infty} \widetilde{\Phi}_{\nu_n, \varepsilon_n}^I(\widetilde{\eta}_n, \check{\theta}_n) = \Lambda_* < \infty
    \end{align*}
    because the other case should be obvious. 
    From \eqref{gamma00} and Proposition \ref{rem-2} (I), we deduce that
    \begin{align}\label{gamma03}
        & \liminf_{n \to \infty} \int_I \int_\Omega \alpha(\widetilde{\eta}_n)|\nabla \check{\theta}_n| \, dx dt \geq \int_I \int_\Omega d \bigl[ \alpha(\widetilde{\eta}(t))|D \check{\theta}(t)| \bigr].
    \end{align}
    We verify the condition of the lower bound
    \begin{align*}
        & \liminf_{n \to \infty} \widetilde{\Phi}_{\nu_n, \varepsilon_n}^I(\widetilde{\eta}_n, \check{\theta}_n) \geq \widetilde{\Phi}_{0}^I(\widetilde{\eta}, \check{\theta}),
    \end{align*}
    as a consequence of \eqref{gamma02}, \eqref{gamma03}, and the lower semi-continuity of the $ L^2 $-based norms. 

    Next, for the condition of optimality, we consider a function $ \hat{\theta} \in D( \widetilde{\Phi}_{0}^I(\widetilde{\eta}, \cdot)) $ $ (\subset L^2(I; X)) $. Then, because Remark \ref{Rem.Note05} (Fact\,4) suggests that $ |D \hat{\theta}(\cdot)|(\Omega) \in L^1(I) $, we can apply the intermediate type approximation methods (cf. Proposition \ref{PropKen01}) and the standard approximation method of the Sobolev function (cf. \cite{MR1857292,MR3409135}) to determine a sequence $ \{ \omega_k \}_{k = 1}^\infty \subset C^\infty(\overline{I \times \Omega}) $ such that
    \begin{subequations}\label{gamma20}
        \begin{align}\label{gamma20a}
            \omega_k \to \hat{\theta} & \mbox{ in $ L^2(I; X) $,} ~\nu_0 \omega_k \to \nu_0 \hat{\theta} \mbox{ in $ L^2(I; Y) $}
        \\
        & \qquad |D \omega_k(\cdot)|(\Omega) \to |D \check{\theta}(\cdot)|(\Omega) \mbox{ in $ L^1(I) $} \mbox{ as $ k \to \infty $.}
            \nonumber
    \end{align}
We consider a sequence $ \{ 1 < \hat{n}_1 < \hat{n}_2 < \hat{n}_3 < \dots < \hat{n}_k < \dots \} \subset \N  $ such that
        \begin{align}\label{gamma20b}\noeqref{gamma20b}
        & \frac{|\nu_n^2 -\nu_0^2|}{2} \bigl| \nabla \omega_k \bigr|_{L^2(I; [X]^N)}^2 \leq 2^{-(k +1)} \mbox{ for any $ n \geq \hat{n}_k $}
    \end{align}
    and define
    \begin{align}\label{gamma20c}
        \hat{\theta}_n := & \left\{ \begin{array}{ll}
            \omega_k, & \mbox{if $ \hat{n}_k \leq n < \hat{n}_{k +1} $ for $ k = 1, 2, 3, \dots $,}
            \\[1ex]
            \omega_1, & \mbox{otherwise}
        \end{array} \right. 
        \\
        & \qquad \mbox{in $ L^2(I; X) $ for $ n = 1, 2, 3, \dots.  $}
        \nonumber
    \end{align}
    \end{subequations}
    From \eqref{gamma20a}--\eqref{gamma20c}, we observe that
    \begin{align}\label{gamma22}
        & \hat{\theta}_n \to \hat{\theta} \mbox{ in $ L^2(I; X) $}, ~ \nu_n \hat{\theta}_n \to \nu_0 \hat{\theta} \mbox{ in $ L^2(I; Y) $},
        \nonumber
        \\
        & \mbox{and } |D \hat{\theta}_n ( \cdot)|(\Omega) \to |D \hat{\theta}(\cdot)|(\Omega) \mbox{ in $ L^1(I) $ as $ n \to \infty $.}
    \end{align}
    Additionally, because of \eqref{gamma00}, we can apply Proposition \ref{rem-2} (II) to the case when
    \begin{align*}
        & \beta = 1, ~ \{ \beta_n \}_{n = 1}^\infty = \{1, 1, 1, \dots\}, ~ \varrho = \alpha(\widetilde{\eta}), ~\mbox{and}~ \{ \varrho_n \}_{n = 1}^\infty = \{ \alpha(\widetilde{\eta}_n) \}_{n = 1}^\infty
    \end{align*}
    and deduce that
    \begin{align}\label{gamma21}
        & \int_I \int_\Omega \alpha(\widetilde{\eta}_n(t)) |\nabla \hat{\theta}_n(t)| \, dx dt \to \int_I \int_\Omega d[\bigl[ \alpha(\widetilde{\eta}(t)) |D \hat{\theta}(t)| \bigr] \, dt \mbox{ as $ n  \to \infty $.}
    \end{align}
    Moreover, bearing \eqref{gamma21} and (\hyperlink{A6l}{A6}) in mind, we compute the following:
    \begin{align}
        & \left| \int_I \int_\Omega \alpha(\widetilde{\eta}_n(t)) \gamma_{\varepsilon_n}(\nabla \hat{\theta}_n(t)) \, dx dt -\int_I \int_\Omega  d[\bigl[ \alpha(\widetilde{\eta}(t)) |D \hat{\theta}(t)| \bigr] \, dt \right| 
        \nonumber
        \\
        & \qquad \leq \left|  \int_I \int_\Omega \alpha(\widetilde{\eta}_n(t)) |\nabla \hat{\theta}_n(t)| \, dx dt -\int_I \int_\Omega d[\bigl[ \alpha(\widetilde{\eta}(t)) |D \hat{\theta}(t)| \bigr] \, dt\right|
        \nonumber
        \\
        & \qquad \qquad +\varepsilon_n \int_I \int_\Omega \alpha(\widetilde{\eta}_n(t)) \, dx dt
        \nonumber
        \\
        & \qquad \to 0 ~\mbox{ as $ n \to \infty $.}
        \label{gamma23}
    \end{align}
    We verify the condition of optimality
    \begin{align*}
        & \widetilde{\Phi}_{\nu_n, \varepsilon_n}^I(\widetilde{\eta}_n, \hat{\theta}_n) \to \widetilde{\Phi}_{0}^I(\widetilde{\eta}, \hat{\theta}) \mbox{ as $ n \to \infty $}
    \end{align*}
    as a consequence of \eqref{gamma22} and \eqref{gamma23}.

    Finally, we prove item (II). Because of \eqref{convOpt} and \eqref{gamma00}, we can apply Proposition \ref{rem-2} (I) to observe that
    \begin{subequations}\label{ken40}
    \begin{align}\label{ken40-01}
        \liminf_{n \to \infty} \int_I \int_\Omega \alpha(\widetilde{\eta}_n) \gamma_{\varepsilon_n}(\nabla \widetilde{\theta}_n) \, dx dt & \geq \liminf_{n \to \infty} \int_I \int_\Omega \alpha(\widetilde{\eta}_n) |\nabla \widetilde{\theta}_n| \, dx dt
        \nonumber
        \\
        & \geq \int_I d \bigl[ \alpha(\widetilde{\eta}(t)) |D \widetilde{\theta}(t)| \bigr] \, dt.
    \end{align}
    Additionally, by the lower semi-continuity of the $L^2$-based norm, it immediately follows that
    \begin{align}\label{ken40-02}
        & \begin{cases}
            \ds \liminf_{n \to \infty} \frac{1}{2} \bigl| \nabla (\nu_n \widetilde{\theta}_n) \bigr|_{L^2(I; [X]^N)}^2 \geq \frac{1}{2} \bigl| \nabla (\nu_0 \widetilde{\theta}) \bigr|_{L^2(I; [X]^N)}^2,
            \\[2ex]
            \ds \liminf_{n \to \infty} \frac{M_0}{2} \bigl| \widetilde{\theta}_n \bigr|_{L^2(I; X)}^2 \geq \frac{M_0}{2} \bigl| \widetilde{\theta} \bigr|_{L^2(I; X)}^2.
        \end{cases}
    \end{align}
    \end{subequations}
    Because of \eqref{convOpt}, \eqref{ken40-01}, and \eqref{ken40-02}, we infer that
    \begin{subequations}\label{ken41}
    \begin{align}
        & \lim_{n \to \infty} \int_I \int_\Omega \alpha(\widetilde{\eta}_n(t)) \gamma_{\varepsilon_n}(\nabla \widetilde{\theta}_n(t)) \, dx dt = \int_I \int_\Omega d \bigl[ \alpha(\widetilde{\eta}(t)) |D \widetilde{\theta}(t)| \bigr] \, dt,
        \label{ken41-01}
        \\[2ex]
        & \lim_{n \to \infty} \frac{1}{2} \bigl| \nabla (\nu_n \widetilde{\theta}_n) \bigr|_{L^2(I; [X]^N)}^2 = \frac{1}{2} \bigl| \nabla (\nu_0 \widetilde{\theta}) \bigr|_{L^2(I; [X]^N)}^2.
        \label{ken41-02}
    \end{align}
    \end{subequations}
        Additionally, from \eqref{ken41-01} and (\hyperlink{A6l}{A6}), we compute the following:
        \begin{align}\label{ken42}
            & \lim_{n \to \infty} \int_I \int_\Omega \alpha(\widetilde{\eta}_n(t))|\nabla \widetilde{\theta}_n(t)| \, dx dt = \lim_{n \to \infty} \int_I \int_\Omega \alpha(\widetilde{\eta}_n(t)) \gamma_{\varepsilon_n}(\nabla \widetilde{\theta}_n(t)) \, dx dt
            \\
            & = \int_I \int_\Omega d \bigl[ \alpha(\widetilde{\eta}(t)) |D \widetilde{\theta}(t)| \bigr] \, dt.
        \end{align}

        Now, we obtain the convergence \eqref{kenConv01} by applying Proposition \ref{rem-2} to the case when
        \begin{align*}
            & \beta = \alpha(\widetilde{\eta}), ~ \{\beta_n\}_{n = 1}^\infty = \{ \alpha(\widetilde{\eta}_n) \}_{n = 1}^\infty, ~ \varrho = 1, ~ \{ \varrho_n \}_{n = 1}^\infty = \{1, 1, 1, \dots\}.
        \end{align*}
        Meanwhile, because \eqref{convOpt} and \eqref{ken41-02} lead to
        \begin{align*}
            & \nu_0 \widetilde{\theta}_n \to \nu_0 \widetilde{\theta} \mbox{ weakly in $ L^2(I; Y) $ as $ n \to \infty $,}
        \end{align*}
        we verify the convergence \eqref{kenConv02} as a consequence of \eqref{ken41-02} and the uniform convexity of the $ L^2 $-based norm. \qed
}
\bigskip

Now, we can prove the Main Theorem. 
\bigskip

\noindent
\textbf{Proof of the Main Theorem. }
We assume (\hyperlink{A5l}{A5}) and (\hyperlink{A7l}{A7}), and for every $ \nu \in (0, \nu_0 +1] $ and $ \varepsilon \in (0, 1) $, we let $ \bigl[ \{ \eta_{\nu, \varepsilon, i}^{m} \}_{i = 1}^{m}, \{ \theta_{\nu, \varepsilon, i}^{m} \}_{i = 1}^{m} \bigr] \subset  [X]^{m \times 2}  $ denote the periodic solution to the approximating problem (S)$_{\nu, \varepsilon}^{m}$, which is obtained in Theorem \ref{Thm.1} with the initial value $ [\eta_{\nu, \varepsilon, 0}^{m}, \theta_{\nu, \varepsilon, 0}^{m}] \in [Y]^2 $.  Then, invoking \eqref{thm1-01} in the proof of Theorem \ref{Thm.1}, we observe that
\begin{description}
    \item[$(*1)$]the class of initial values $ \left\{ \begin{array}{l|l} 
            [\eta_{\nu, \varepsilon, 0}^{m}, \theta_{\nu, \varepsilon, 0}^{m}]  & m \in \mathbb{N}, \nu \in (0, \nu_0 +1], \varepsilon \in (0, 1)
    \end{array} \right\} $  is contained in a compact convex set $ \mathbb{K}_0 $ in $ [X]^2 $ defined as
\begin{align*}
    \mathbb{K}_0 ~& := \left\{ \begin{array}{l|l}
        [\tilde{\eta}, \tilde{\theta}] \in Y \times BV(\Omega) & 
         \parbox{5.5cm}{
             $ |\tilde{\eta}|_{L^\infty(\Omega)} \vee |\tilde{\theta}|_{L^\infty(\Omega)} \leq R_0 $, 
             \\[1ex]
             $ \ds |\tilde{\eta}|_X^2 + R_*|\tilde{\theta}|_X^2 \leq R_1 $, 
             \\[1ex] 
            $ \kappa^2 |\nabla \tilde{\eta}|_{[X]^N}^2 +\delta_*|D\tilde{\theta}|(\Omega) \leq  R_3 $
    }
    \end{array} \right\}.
\end{align*}
\end{description}
Hence, there exist sequences $\{ m_n\}_{n=1}^\infty \subset \mathbb{N} $, $ \{ \nu_n \}_{n = 1}^\infty \subset (0, \nu_0 +1] $ and $ \{ \varepsilon_n \}_{n = 1}^\infty \subset (0, 1) $, and a pair of functions $ [\eta_0, \theta_0] \in \mathbb{K}_0 $ such that
\begin{align*}
    & m_n \uparrow \infty ~\mbox{as $ n \to \infty $,} ~\mbox{and}~ \{\nu_n\}_{n=1}^\infty \mbox{ and } \{\varepsilon_n\}_{n=1}^\infty \mbox{ satisfy } \eqref{bar-nu-eps01}, 
\end{align*}
and
\begin{align*}
    & \begin{cases}
        \eta_{n, 0} := \eta_{\nu_n, \varepsilon_n, 0}^{m_n} \to \eta_0 ~\mbox{in $X $, weakly in $ Y $,}
        \\[1ex]
        \theta_{n, 0} := \theta_{\nu_n, \varepsilon_n, 0}^{m_n} \to \theta_0 ~\mbox{in $ X $, weakly-$*$ in $ BV(\Omega) $}
    \end{cases}
    \mbox{as $ n \to \infty $.}
\end{align*}

Next, 
we construct the following sequences of different time interpolations:
\begin{align*}
    & \begin{cases}
        \bigl[ \overline{\eta}_n(t), \overline{\theta}_n(t) \bigr] := [\eta_{\nu_n, \varepsilon_n, i}^{m_n}, \theta_{\nu_n, \varepsilon_n, i}^{m_n}], ~ \bigl[ \underline{\eta}_n(t), \underline{\theta}_n(t) \bigr] := [\eta_{\nu_n, \varepsilon_n, i-1}^{m_n}, \theta_{\nu_n, \varepsilon_n, i-1}^{m_n}],
        \\[1ex]
        \ds \bigl[ {\eta}_n(t), {\theta}_n(t) \bigr] := \frac{i\tau_n -t}{\tau_n} [\eta_{\nu_n, \varepsilon_n, i-1}^{m_n}, \theta_{\nu_n, \varepsilon_n, i-1}^{m_n}] +\frac{t -(i -1)\tau_n}{\tau_n} [\eta_{\nu_n, \varepsilon_n, i}^{m_n}, \theta_{\nu_n, \varepsilon_n, i}^{m_n}],
    \end{cases}
    \mbox{in $ [X]^2 $}
\end{align*}
and
\begin{align}\label{force00}
    & \bigl[ \overline{u}_n(t), \overline{v}_n(t) \bigr] := [u_{i}^{m_n}, v_{i}^{m_n}] ~\mbox{in $ [X]^2 $}
\end{align}
\begin{center}
for all $ t \in [(i -1)\tau_n, i \tau_n) $, $ i = 1, \dots, m_n $, and  $ n = 1, 2, 3, \dots $.
\end{center}

We fix the arbitrary open interval $ I \subset (0, T) $. Then, from Definition \ref{Def.time-dis} and Remark \ref{Shorem01}, we observe that
\begin{align}\label{tInterp_a}
    & \begin{cases}
        \bigl\{ [\overline{\eta}_n, \overline{\theta}_n],\, [\underline{\eta}_n, \underline{\theta}_n] \bigr\}_{n = 1}^\infty \subset \bigl[ L^\infty(I; Y) \bigr]^2,
        \\[1ex]
        \bigl\{ [{\eta}_n, {\theta}_n] \bigr\}_{n = 1}^\infty \subset \bigl[ W^{1, 2}(I; Y) \bigr]^2,
    \end{cases}
\end{align}
\begin{align}\label{tInterp_c}
    \int_I \bigl( \partial_t \eta_n & (t), \varphi(t) \bigr)_X \, dt +\kappa^2 \int_I \bigl( \nabla \overline{\eta}_n(t), \nabla \varphi(t) \bigr)_{[X]^N} \, dt +\int_I \bigl( g(\overline{\eta}_n(t)), \varphi(t) \bigr)_X \, dt
    \nonumber
    \\
    & +\int_I \varphi(t) \alpha'(\overline{\eta}_n(t)) \gamma_{\varepsilon_n}(\nabla \overline{\theta}_n(t)) \, dt = 0 ~\mbox{for any $ \varphi \in L^2(I; Y) $,}
\end{align}
\begin{subequations}\label{tInterp_bd}
\begin{align}\label{tInterp_b}\noeqref{tInterp_b}
    \int_I \bigl( \alpha_0(\underline{\eta}_n(t))\,&  \partial_t \theta_n(t) -\overline{v}_n(t), \overline{\theta}_n(t) -\psi(t) \bigr)_X \, dt +\widetilde{\Phi}_{\nu_n, \varepsilon_n}^{I}(\underline{\eta}_n; \overline{\theta}_n) 
    \nonumber
    \\
    & \leq \widetilde{\Phi}_{\nu_n, \varepsilon_n}^{I}(\underline{\eta}_n; \psi) ~\mbox{for any $ \psi \in L^2(I; Y) $;}
\end{align}
\begin{align}\label{tInterp_d}\noeqref{tInterp_d}
    \mbox{that is,} ~~& \bigl[ \, \overline{\theta}_n, ~\alpha_0(\underline{\eta}_n) \partial_t \theta_n -\overline{v}_n  \bigr] \in \partial \widetilde{\Phi}_{\nu_n, \varepsilon_n}^I(\overline{\eta}_n; \cdot ) ~\mbox{in $ L^2(I; X) \times L^2(I; X) $.}
\end{align}
\end{subequations}
Additionally, as a consequence of Lemma \ref{Lem3-01}, we infer that
\begin{description}
    \item[$(*2)$]the sequences $ \bigl\{ \overline{\eta}_n \bigr\}_{n = 1}^\infty $ and $ \bigl\{ \underline{\eta}_n \bigr\}_{n = 1}^\infty $ are bounded in $ L^\infty(I; Y) \cap L^\infty(Q) $, and the sequence $ \bigl\{ \eta_n \bigr\}_{n = 1}^\infty $ is bounded in $ W^{1, 2}(0, T; X) \cap L^\infty(I; Y) \cap L^\infty(Q) $;
    \item[$(*3)$]the sequences $ \bigl\{ \overline{\theta}_n \bigr\}_{n = 1}^\infty $ and $ \bigl\{ \underline{\theta}_n \bigr\}_{n = 1}^\infty $ are bounded in $ L^\infty(I; X) \cap L^\infty(Q) $, the sequence $ \bigl\{ \eta_n \bigr\}_{n = 1}^\infty $ is bounded in $ W^{1, 2}(0, T; X) \cap L^\infty(Q) $, and the sequences $ \bigl\{ |D \overline{\theta}_n(\cdot)|(\Omega) \bigr\}_{n = 1}^\infty $, $ \bigl\{ |D \underline{\theta}_n(\cdot)|(\Omega) \bigr\}_{n = 1}^\infty $, and $ \bigl\{ |D {\theta}_n(\cdot)|(\Omega) \bigr\}_{n = 1}^\infty $ are bounded in $ L^\infty(I) $.
\end{description}
Therefore, bearing the compactness theories in \cite[Theorem 5.7.7]{MR2192832} and \cite[Corollary 4]{MR0916688} in mind, we can assume the existence of subsequences of $ \bigl\{ [\overline{\eta}_n, \overline{\theta}_n] \bigr\}_{n = 1}^\infty $, $ \bigl\{ [\underline{\eta}_n, \underline{\theta}_n] \bigr\}_{n = 1}^\infty $, and $ \bigl\{ [{\eta}_n, {\theta}_n] \bigr\}_{n = 1}^\infty $ (not relabeled), together with a limiting pair $ [\eta, \theta] \in [L^2(I; X)]^2 $, such that
\begin{subequations}\label{limits}
\begin{align}\label{limits_a}
    & \begin{cases}
        \hspace{-2.5ex}\parbox{12cm}{
            \vspace{-1.5ex}
            \begin{itemize}
                \item $ \overline{\eta}_n \to \eta $ and $ \underline{\eta}_n \to \eta $ in $ L^\infty(I; X) $ and weakly-$*$ in $ L^\infty(I; Y) $,
                    \vspace{-1ex}
                \item $ {\eta}_n \to \eta $  in $ C(\overline{I}; X) $, weakly in $ W^{1, 2}(I; X) $, and weakly-$*$ \linebreak in $ L^\infty(I; Y) $,
            \vspace{-1.5ex}
            \end{itemize}
    }
    \end{cases}
\end{align}
    \vspace{-1ex}
\begin{align}\label{limits_b}
    & \begin{cases}
        \hspace{-2.5ex}\parbox{12cm}{
            \vspace{-1.5ex}
            \begin{itemize}
                \item $ \overline{\theta}_n \to \eta $ and $ \underline{\theta}_n \to \eta $ in $ L^\infty(I; X) $,
                    \vspace{-1ex}
                \item $ {\theta}_n \to \theta $  in $ C(\overline{I}; X) $, and weakly in $ W^{1, 2}(I; X) $
            \vspace{-1.5ex}
            \end{itemize}
    }
    \end{cases}
\end{align}
as $ n \to \infty $. In particular, because of the periodicity of the approximating solutions,
\begin{align}\label{limits_c}
    [\eta(0), \theta(0)] ~& = \lim_{n \to \infty} [\eta_{n, 0}, \theta_{n, 0}] = \lim_{n \to \infty} [\eta_{n}^{m_n}, \theta_{n}^{m_n}] 
    \nonumber
    \\
    & = [\eta(T), \theta(T)] ~\mbox{ in $ [X]^2 $.}
\end{align}
\end{subequations}
Furthermore, bearing \eqref{force00}, (\hyperlink{A2l}{A2}), (\hyperlink{A3l}{A3}), and (\hyperlink{A7l}{A7}) in mind, we observe that
\begin{align}\label{force01}
    & -\alpha_0(\underline{\eta}_n) \partial_t \theta_n +\overline{v}_n \to  -\alpha_0(\eta) \partial_t \theta +v ~\mbox{weakly in $ L^2(I; X) $ as $ n \to \infty $.}
\end{align}

Now, in \eqref{tInterp_bd}, we take the limit as $ n \to \infty $. Then, from Lemma \ref{Lem3-Gamma} and Remark \ref{Rem.M-Gconvs}, we deduce that
\begin{align}\label{fin01}
    & \bigl[ \theta, \, -\alpha_0(\eta) \partial_t \theta +v  \bigr] \in \partial {\widetilde{\Phi}_{0}}^I(\eta; \cdot ) ~\mbox{in $ L^2(I; X) \times L^2(I; X) $}
\end{align}
and
\begin{align}\label{fin02}
    & \int_I |D \overline{\theta}_n(t)|(\Omega) \, dt \to \int_I |D \theta(t)|(\Omega) \, dt ~\mbox{as $ n \to \infty $. }
\end{align}
Because the choice of the open interval $ I \subset (0, T) $ is arbitrary, 
\begin{description}
    \item[$(*4)$]the inclusion of \eqref{fin01} and Remark \ref{Rem.2ndEq} lead to the variational inequality \eqref{S3}.
\end{description}

In the meantime, considering arbitrary $ \varphi \in L^2(I; Y) $ and applying Proposition \ref{rem-2} to the case when\begin{align*}
        & \beta = 1, ~ \{ \beta_n \}_{n = 1}^\infty = \{1, 1, 1, \dots \}, ~ \varrho = \varphi \alpha(\eta), ~\mbox{and}~ \{ \varrho_n \}_{n = 1}^\infty = \{ \varphi \alpha(\overline{\eta}_n) \}_{n = 1}^\infty,
\end{align*}
we observe from \eqref{limits_a}, \eqref{limits_b}, \eqref{fin02}, and (\hyperlink{A3l}{A3}) that
\begin{align}\label{fin03}
    & \left| \int_I \int_\Omega \varphi(t) \alpha(\overline{\eta}_n(t)) \gamma_{\varepsilon_n}(\nabla \overline{\theta}_n(t)) \, dx dt -\int_I \int_\Omega d \bigl[ \varphi(t) \alpha({\eta}(t)) |\nabla {\theta}(t)| \bigr] \, dt \right|
    \nonumber
    \\[2ex]
    & \qquad \leq \left| \int_I \int_\Omega \varphi(t) \alpha(\overline{\eta}_n(t)) |\nabla \overline{\theta}_n(t)| \, dx dt -\int_I \int_\Omega d \bigl[ \varphi(t) \alpha({\eta}(t)) |\nabla {\theta}(t)| \bigr] \, dt \right|
    \nonumber
    \\[1ex]
    & \qquad \qquad +\varepsilon_n \sup_{n \in \N} \bigl| \varphi \alpha(\overline{\eta}_n) \bigr|_{L^1(Q)}
    \nonumber
    \\[2ex]
    & \qquad \to 0 ~\mbox{ as $ n \to \infty $.}
\end{align}
Because of \eqref{limits_a}, \eqref{limits_b}, \eqref{fin03}, and (\hyperlink{A2l}{A2}), letting $ n \to \infty $ in \eqref{tInterp_c} yields \begin{align}\label{fin04}
    \int_I \bigl( \partial_t \eta & (t), \varphi(t) \bigr)_X \, dt +\kappa^2 \int_I \bigl( \nabla {\eta}(t), \nabla \varphi(t) \bigr)_{[X]^N} \, dt +\int_I \bigl( g({\eta}(t)), \varphi(t) \bigr)_X \, dt
    \nonumber
    \\
    & \hspace{10ex}+\int_I \int_\Omega d \bigl[\varphi(t) \alpha'({\eta}(t)) |D \overline{\theta}(t)| \bigr] \, dt = 0 ~
    \\[1ex]
    & \hspace{3ex}\mbox{for any $ \varphi \in L^2(I; Y) $, and any open interval $ I \subset (0, T) $. }
    \nonumber
\end{align}
\begin{description}
    \item[$(*4)$]\eqref{fin04} is equivalent to the variational identity \eqref{S2}.
\end{description}

$(*4)$, $(*5)$, and \eqref{limits_c} complete the proof of the Main Theorem. 
\qed

\paragraph*{Acknowledgements}{
This work was supported by Grant-in-Aid for Scientific Research (C) No. 20K03672, JSPS.  
We thank Edanz (https://jp.edanz.com/ac) for editing a draft of this manuscript.
}



\end{document}